%% file: preprint.tex
\documentclass[a4paper]{article}

\usepackage{mystyle}
\usepackage{review}
\usepackage{esint}
\usepackage{standalone}

\title{Combined Regularization and Discretization of Equilibrium Problems and Primal-Dual Gap Estimators}

\author{Steven-Marian Stengl}

\newcommand\my{\operatorname{MY}}
\newcommand\bR{\overline{\mathbb{R}}}
\newcommand\Lip[1]{\operatorname{Lip}\left(#1\right)}
\newcommand\lowobs{\underline{\psi}}
\newcommand\uppobs{\overline{\psi}}

\newcommand\moscoconv{\ensuremath{\overset{M}{\longrightarrow}}}
\newcommand\lowobsh{\ensuremath{\underline{\psi}_\calT}}
\begin{document}
\maketitle
\abstract{The present work aims at the \msa{application}{} of finite element discretizations \msa{to}{} a class of equilibrium problems involving moving constraints. Therefore, a Moreau--Yosida based regularization technique, controlled by a parameter, is discussed \msa{and, using a}{} generalized $\Gamma$-convergence concept\msa{,}{} a priori \msa{convergence}{} results \msa{are derived}{}. The \msa{latter}{} technique is applied to the discretization \msa{of the regularized problems and is used to prove the convergence to the orginal equilibrium problem, when}{} both \msa{--- regularization and discretization --- are imposed simultaneously}. In addition, a primal-dual gap technique is used for the derivation of error estimators \msa{suitable for adaptive mesh refinement}. A strategy for balancing between a refinement of the mesh and an update of the regularization parameter is established\msa{, too}. The theoretical findings are illustrated for the obstacle problem as well as numerical experiments are performed for two quasi-variational inequalities with application to thermoforming and biomedicine, respectively.}
\section{Introduction}
In analysis, a broad class of problems such as optimization problems, Nash games, variational inequalities (VI) as well as their more general counterpart quasi-variational inequalities (QVI) can be enveloped under the umbrella of \emph{(quasi-)equilibrium problems} (see \cite{bib:BlumOettliEquilibrium, bib:AusselCotrinaIusemQuasiEquilibrium}). Especially (Q)VIs, originally introduced in \cite{bib:BensoussanLionsImpulseFrench} proved themselves to be a powerful tool for modelling a variety of applications ranging from superconductivity (cf. \cite{bib:BarrettPrigozhinSuperconductor}) over sandpile formation (cf. \cite{bib:BarrettPrigozhinSandpile}) to technical processes like thermoforming (cf. \cite{bib:AlphonseHintRautenbergDirectDiff, bib:AlphonseRautenbergRodrigues}). All these models live in a function space setting. \msa{As the derivation of exact solutions is often beyond reach}{}, discretization methods such as finite element methods \msa{are}{} applied for the \msa{numerical approximation}{} of solutions. This technique has been successfully used for a selection of VIs (cf. \cite{bib:FalkErrorEstimateVI, bib:BrezziHagerRaviartPartOne, bib:BrezziHagerRaviartPartTwo}) as well as QVIs (cf. \cite{bib:BarretPrigozhinSandpilesNonconformingFEM}). Especially for the latter case the literature on adaptive refinement of the underlying mesh is scarce (cf. \cite{bib:BartelsCarstensenAveragingAPosterioriVI, bib:BraessAPosterioriObstacle} as well as the overview article \cite{bib:WohlmuthOverview} and the references therein for variational inequalities).
The present work addresses this and discusses in an abstract framework for equilibrium problems a priori convergence results and derives error estimators based on the techniques in \cite{bib:BartelsMilicevicPrimalDualGap}.\\
\msa{To obtain a clearer view of the problems addressed in the scope of this work we draw our attention to the following \emph{obstacle problem} (cf. \cite{bib:KinderlehrerStampacchia,bib:Rodrigues}):}\\
For a given forcing term $f \in L^2(\Omega)$ and an obstacle function $\lowobs \in H^1(\Omega)$ with $\lowobs|_{\partial \Omega} < 0$ a.e. on $\partial \Omega$ we \msa{define}{} the constraint set
\begin{eq*}
K:= \{z \in H_0^1(\Omega) : \msa{z}{} \geq \lowobs \text{ a.e. on } \Omega\}.
\end{eq*}
The\msa{n, the}{} obstacle problem reads as
\begin{eq}\label{eq:obstacleproblem}
\minimize\hspace*{-0.7em} \frac{1}{2}\|\nabla y\|_{L^2(\Omega; \R^d)}^2 - (f,y)_{L^2(\Omega)} \, \subjectto y \in K.
\end{eq}
\msa{The constraint condition can equivalently be rewritten by adding the \emph{indicator function} of $K$ defined by
\begin{eq*}
I_K(y) := \left\{\begin{array}{ll} 0, &\text{if } y \in K,\\ \infty, & \text{else.}  \end{array}\right.
\end{eq*}}
\msa{The problem \refer{eq:obstacleproblem} is frequently studied in non-smooth and convex analysis as it includes aspects relevant for other applications, such as the presence of a constraint leading to multipliers of low regularity. Moreover, the first order system reads as the following VI:\\
Seek $y \in K$ such that
\begin{eq}[\label{eq:obstacleproblem:vi}]
(\nabla y, \nabla v - \nabla y)_{L^2(\Omega;\R^d)} \geq (f,v-y)_{L^2(\Omega)} \text{ for all } v \in K.
\end{eq}
For the connection of \refer{eq:obstacleproblem} and \refer{eq:obstacleproblem:vi} the interested reader is referred to \cite[Section 1:2]{bib:Rodrigues}. The obstacle problem serves as a prototypical VI as well as an optimization problem.}{}
As an \msa{instance of the latter}{} it falls as well in the category of \msa{\emph{equilibrium problems}}.
\msa{Introducing the latter we refer to \cite{bib:HintStenglEquiGamma} and the references therein.}{}\\
\msa{By $\bR$ we denote the set $\R \cup \{+\infty\}$.}{} Let $U$ be a reflexive Banach space and let $U_\ad$ denote a non-empty, closed, convex subset of $U$. Let a functional $\calE: U_\ad \times U_\ad \rightarrow \overline{\mathbb{R}}$ with  $\dom{\calE(\ccdot,u)} \neq \emptyset$ for all $u \in U_\ad$ be given. A point $u \in U_\ad$ is called \emph{equilibrium}  (see \msa{\cite[Definition 1]{bib:HintStenglEquiGamma}}), if
\begin{eq}\label{eq:equi}
\calE(u,u) \leq \calE(v,u) \text{ holds for all } v \in U_\ad.
\end{eq}
In the scope of this work we might refer to the first component as \emph{control component} and \msa{to}{} the second component as \emph{feedback component}. 
This concept is a slight deviation from the \msa{ones}{} in \cite{bib:BlumOettliEquilibrium, bib:AusselCotrinaIusemQuasiEquilibrium}. However, in \msa{here}{} we are actually interested in minimizers of the functional $\calE(\ccdot,u)$. \msa{The type of functionals investigated in this work takes the form}{}
\begin{eq}\label{eq:intro}
\calE(v,u) := F(v,u) + G(Av,u) + H(Bv,u)
\end{eq}
with functionals $F:U \times U \to \bR$, $G: Y \times U \to \bR$ and $H:Z \times U \to \bR$. Here, $Y$ is another reflexive Banach space and $Z$ a (real) Hilbert space and $A \in \calL(U,Y), B \in \calL(U,Z)$ are bounded, linear operators on their \msa{respective}{} spaces. Of course, the obstacle problem fits into that framework, which we will elaborate in more detail in the main body of \msa{this}{} work. However, the presence of  the obstacle constraint poses a challenge. In \msa{this work}{}, the presence of a non-smooth or even non-continuous part in the functional, such as the presence of the constraint\msa{,}{} will be represented \msa{by}{} the functional $H$.\\ 
One popular technique in use is the substitution of the problem with a sequence of more easily treatable problems. For this sake a suitable convergence concept needs to be applied. In optimization, the concept of $\Gamma$-convergence proved to be a versatile tool especially for phase field approximations \cite{bib:BraidesGammaConvBeginners}. Here, we give the following extension (cf. \cite[Definitions 8 and 9]{bib:HintStenglEquiGamma}):\\
Given \msa{functionals}{} $\msa{(}\calE_n\msa{)_{n \in \N}, \calE}:U_\ad \times U_\ad \to \bR$, we say that $\msa{(}\calE_n\msa{)_{n \in \N}}{}$ is (weakly) \emph{$\Gamma$-convergent} \msa{to}{} $\calE$, if for all sequences $u_n \to u$ (resp. $u_n \rightharpoonup u$) holds
\begin{eq*}
\calE(u,u) \leq \liminf_{n \to \infty}\calE_n(u_n,u_n)
\end{eq*}
and for all $v \in U_\ad$ there exists $v_n \to v$ (resp. $v_n \rightharpoonup v$) such that
\begin{eq*}
\limsup_{n \to \infty} \calE_n(v_n,u_n) \leq \calE(v,u).
\end{eq*}
This generalizes the notion of $\Gamma$-convergence to equilibrium problems and is closely related to the corresponding concept \msa{\emph{multi-epiconvergence}}{} in \cite{bib:JSPangGuerkan} that \msa{has}{} been developed for (generalized) Nash equilibrium problems. \msa{The}{} notion of Mosco-convergence is generalized in a similar way\msa{:}{}\\
A sequence $(\calE_n)_{n \in \N}$ is called \emph{Mosco-convergent}, if the first condition with respect to the weak convergence holds and for all $u_n \rightharpoonup u$ and $v \in U_\ad$ there exists a sequence $v_n \to v$ such that
\begin{eq*}
\limsup_{n \to \infty} \calE_n(v_n,u_n) \leq \calE(v,u).
\end{eq*}
\msa{For optimization problems}{} $\Gamma$-convergence enjoys numerical interest, since accumulation points of minimizers of the iterates are as well minimizers of the limit. An analogous result can be cited for equilibrium problems as well in the following \msa{result}.
\begin{thm}[see \mbox{\msa{\cite[Theorem 11]{bib:HintStenglEquiGamma}}}]\label{prop:equi:gamma:mini_convergence}
Let $(\calE_n)_{n \in \N}$ be a (weakly) $\Gamma$-con\-ver\-gent sequence of functionals with limit \msa{$\calE$}. Then, every (weak) accumulation point of a sequence of corresponding equilibria $(u_n)_{n \in \N}$ is an equilibrium of the limit.  
\end{thm}
The rest of this paper is organized as follows: In Section 2 we discuss the employed Moreau--Yosida regularization strategy as well as its $\Gamma$-convergence. This technique goes along with the simultaneous \msa{presence}{} of a discretization and $\Gamma$-convergence is established \msa{in this case}{} as well. Section 3 addresses the establishment of the formulation of primal dual gap error estimators for equilibrium functionals. Section 4 is then devoted to the application of the developed results and techniques to a selection of obstacle-type quasi-variational inequalities. 
\section{A Priori Convergence for Regularized and Discretized Equilibrium Problems}
As announced in the introduction, we are interested in \msa{equilibrium}{} problems involving pointwise constraints. A successfully applied approach to address these constraints are penalization and regularization techniques \msa{(see e.g.}{} \cite{bib:HintItoKunisch,bib:HintKunischPath}\msa{)}. A popular \msa{instance}{} in use is the \emph{Moreau--Yosida regularization}\msa{,}{} whose definition is given next.
\begin{defn}[Moreau--Yosida regularization (cf. \mbox{\cite[Definition 12.20]{bib:BauschkeCombettes}})]\label{defn:my}
Let a (real) Hilbert space $Z$ together with \msa{a}{} positive real number $\gamma > 0$ and a convex, proper, \lsc{} functional $H : Z \to \bR$ be given. The \emph{Moreau--Yosida regularization} of $H$ with respect to $\gamma$ is defined by
\begin{eq*}
\my(\gamma, H)(z) := \inf_{\zeta \in Z} \left( H(\zeta) + \frac{\gamma}{2}\|\zeta - z\|_Z^2 \right).
\end{eq*} 
\end{defn}
As it can be seen from its definition, the Moreau--Yosida regularization is an
infimal convolution (see \cite[Chapter 12]{bib:BauschkeCombettes}) of the
functional $H$ with \msa{the functional}{} $\frac{\msa{\gamma}}{2}\|\cdot\|^2$. It is worth noting, that the
Moreau--Yosida regularization has as its domain always the whole space. The
positive number $\gamma$ serves as \emph{regularization parameter} and controls the influence of the quadratic part. In the light of problems with obstacle-type constraints the functional $H$ might read as an indicator functional of some non-empty, closed, convex subset, with the set possibly depending on the feedback component\msa{, too}. \msa{Returning}{} to \refer{eq:obstacleproblem} one way is the choice $Z = H_0^1(\Omega)$ along with $H := I_K$ and $B = \mathrm{id}_{H_0^1(\Omega)}$. This approach however leads to practical problems, since the calculation of the \msa{Moreau}--Yosida regularization is just as challenging as the obstacle problem itself. Another approach (see \cite{bib:HintKunischPath}) is to use the embedding $H_0^1(\Omega) \hookrightarrow L^2(\Omega)$ and define the set
\begin{eq*}
K_{L^2} := \{ z \in L^2(\Omega) : z \geq \lowobs \text{ a.e. on } \Omega\}
\end{eq*}
and observe, that for $y \in H_0^1(\Omega)$ \msb{the condition $y \in K$ is equivalent to}{} $y \in K_{L^2}$. Hence, we rewrite $I_K = I_{K_{L^2}} \circ B$ with $B = i_{L^2(\Omega)}$  as well as $Z = L^2(\Omega)$ and set $H = I_{K_{L^2}}$ instead.
This allows an explicit calculation of the Moreau--Yosida regularization leading to
\begin{eq*}
\my(\gamma,H)(z) = \frac{\gamma}{2} \int_\Omega (\lowobs - z)^{2+} \dx.
\end{eq*}
%
\subsection{\texorpdfstring{$\Gamma$-Convergence for Moreau--Yosida Regularization of Equilibrium Problems}{Gamma-Convergence for Moreau--Yosida Regularization of Equilibrium Problems}}\label{sec:regularizedgamma}
The general idea is the formulation of a sequence of approximate equilibrium problems via the substitution of the functional $H$ in \refer{eq:intro} by its Moreau--Yosida regularization. For a sequence $(\gamma_n)_{n \in \N} \nearrow \infty$ one expects to recover the original problem. As \msa{convergence notion}{} we utilize the $\Gamma$-conver\-gence as introduced before. This is covered in the following Theorem.
\begin{thm}\label{thm:myreg:gammaconv}
Let functionals $F:U_\ad \times U_\ad \to \bR, G: Y \times U_\ad \to \bR$ and $H:Z \times U \to \bR$\msa{, being  proper, convex and lsc. in the control component,}{} be given. Let for all sequences $u_n \rightharpoonup u$ and $v_n \to v$ hold
\begin{eq*}
F(u,u) &\leq \liminf_{n \to \infty} F(u_n,u_n) &&\text{and}& \limsup_{n \to \infty} F(v_n,u_n) &\leq F(v,u)
\end{eq*}
as well as
\begin{eq*}
G(Au,u) &\leq \liminf_{n \to \infty} G(A u_n,u_n) &&\text{and}& \limsup_{n \to \infty} G(Av_n,u_n) &\leq G(Av,u)
\end{eq*}
and let the assumptions from \refer{lem:my:mosco} hold. Moreover assume, for all sequences $u_n \rightharpoonup u$ and all $v \in U$ there exists a sequence $v_n \rightharpoonup v$ with $\limsup_{n \to \infty} H(Bv_n,u_n) \leq H(Bv,u)$.\\
Define the functionals
\begin{eq*}
\calE(v,u) &:= F(v,u) + G(Av,u) + H(Bv,u) \text{ and}\\
\calE_n(v,u) &:= F(v,u) + G(Av,u) + \my(\gamma_n, H(\ccdot, u))(Bv).
\end{eq*}
Then, the convergence $\calE_n \overset{\Gamma}{\rightharpoonup} \calE$ holds true.
\end{thm}
Before proving this result, we propose the following auxiliary lemmas solely aiming at the behaviour of the Moreau--Yosida regularization itself. As simplest \msa{situation}, we consider the case \msa{when}{} no feedback component occurs.
\begin{lem}\label{lem:my:gamma}
Let a functional $H: Z \to \bR$ be given, then for every sequence $\gamma_n \to \infty$ holds
\begin{eq*}
\my(\gamma_n, H) \moscoconv H.
\end{eq*}
Moreover, for every $z$ the constant recovery sequence $z_n = z$ can be taken.
\end{lem}
\begin{proof}
We utilize the results in \cite[Theorem 3.26]{bib:AttouchEpi}. Therefore, the assertion holds\msa{,}{} if and only if for all $\lambda > 0$ and $z \in Z$ holds
\begin{eq*}
\my(\lambda, \my(\gamma_n, H))(z) \to \my(\lambda,H)(z).
\end{eq*}
Using \cite[Prop. 2.68]{bib:AttouchEpi} we obtain in our notation 
\begin{eq*}
\my(\lambda,\my(\gamma_n,H)) = \my\left(\frac{\lambda \gamma_n}{\lambda + \gamma_n}, H\right) = \my\left(\lambda_n,H\right)
\end{eq*}
with $\lambda_n := \lambda \left(1+ \frac{\lambda}{\gamma_n}\right)^{-1}$. Clearly $\lambda_n \to \lambda$ as $n \to \infty$.\\
For every given $z\in Z$ we see, that the mapping $\lambda' \mapsto \my(\lambda',H)(z) = \inf_{y \in Z}\left( H(y) + \frac{\lambda'}{2}\|y-z\|^2 \right)$ is an infimum of functionals, that are concave in $\lambda'$\msa{. Hence, it is concave as well and as a function on the real numbers also continuous (cf. \cite[Corollary 2.3]{bib:EkelandTemam}).}{} Therefore, \msa{we deduce}{} the pointwise convergence
\begin{eq*}
\my\left(\lambda_n,H\right)(z) \to \my(\lambda,H)(z)
\end{eq*}
and using \cite[Theorem 3.26]{bib:AttouchEpi} the assertion. Moreover, using the inequality $\my(\lambda_n,H)(z) \leq H(z)$ we deduce
\begin{eq*}
\limsup_{n \to \infty} \my(\lambda_n,H)(z) \leq H(z)
\end{eq*}
and hence obtain as possible recovery sequence $z_n = z$. 
\end{proof}
After establishing the $\Gamma$-convergence result for a functional\msa{,}{} solely depending on the control component, we want to turn our attention to the situation involving a feedback component \msa{next}. For this purpose, we propose the following \msa{auxiliary}{} lemma \msa{first}.
\begin{lem}\label{lem:my:conv}
Let a functional $H:Z\times U\to \bR$\msa{, being bounded from below,}{} be given, such that the following two conditions hold\msa{:}{}
\begin{enum}
\item For all sequences $z_n \rightharpoonup z$ and $u_n \rightharpoonup u$ holds
\begin{eq*}
H(z,u) \leq\liminf_{n \to \infty} H(z_n,u_n).
\end{eq*}
\item For all sequences $u_n \rightharpoonup u$ and $z\in Z$ exists a sequence $z_n \to z$ such that
\begin{eq*}
\limsup_{n \to \infty} H(z_n,u_n) \leq H(z,u)
\end{eq*}
holds.
\end{enum}
Then, for every sequence $\lambda_n \to \lambda$ and $z_n \rightharpoonup z$ holds
\begin{eq*}
\my(\lambda,H(\ccdot,u))(z) \leq \liminf_{n \to \infty} \my(\lambda_n,H(\ccdot,u_n))(z_n)
\end{eq*}
and for all $\lambda_n \to \lambda$ and $z_n \to z$ holds
\begin{eq*}
\limsup_{n \to \infty} \my(\lambda_n, H(\ccdot,u_n))(z_n) \leq \my(\lambda,H(\ccdot,u))(z).
\end{eq*}
\end{lem}
The proof is technical and space-consuming and is therefore given in the appendix. However, with \msa{these preparations}{} at hand we are set up to formulate and proof the announced $\Gamma$-convergence result for Moreau--Yosida regularizations in the presence of a feedback component.
\begin{lem}\label{lem:my:mosco}
Let a functional $H:Z\times U\to \bR$\msa{,}{} being bounded from below \msa{be given}{}, such that the following two conditions hold\msa{:}{}
\begin{enum}
\item For all sequences $z_n \rightharpoonup z$ and $u_n \rightharpoonup u$ holds
\begin{eq*}
H(z,u) \leq\liminf_{n \to \infty} H(z_n,u_n).
\end{eq*}
\item For all sequences $u_n \rightharpoonup u$ and $z\in Z$ exists a sequence $z_n \to z$ such that holds
\begin{eq*}
\limsup_{n \to \infty} H(z_n,u_n) \leq H(z,u).
\end{eq*}
\end{enum}
Then, for all sequences $\gamma_n \to \infty$ and $u_n \rightharpoonup u$ the Mosco-convergence 
\begin{eq*}
\my(\gamma_n,H(\ccdot,u_n)) \moscoconv H(\ccdot,u)
\end{eq*}
holds true.
\end{lem}
\begin{proof}
Let arbitrary sequences $u_n \rightharpoonup u$ and $\gamma_n \to \infty$ be given. According to \cite[Theorem 3.26]{bib:AttouchEpi} the assertion is equivalent to the pointwise convergence
\begin{eq*}
\my(\lambda, \my(\gamma_n,H(\ccdot,u_n))) \to \my(\lambda,H(\ccdot,u))
\end{eq*}
for all $\lambda > 0$. Using \cite[Prop. 2.68]{bib:AttouchEpi} we obtain \msa{as in \refer{lem:my:gamma} in}{} our notation
\begin{eq*}
\my(\lambda, \my(\gamma_n,H(\ccdot,u_n))) = \my(\lambda_n , H(\ccdot, u_n))
\end{eq*}
with $\lambda_n = \lambda(1+\frac{\lambda}{\gamma_n})^{-1}$. Clearly, holds $\lambda_n \to \lambda$ as $\gamma_n \to \infty$. Hence, we use \refer{lem:my:conv} and obtain with $z_n = z$ for all $n \in \N$ the convergence
\begin{eq*}
\my(\lambda, H(\ccdot,u))(z) &\leq \liminf_{n \to \infty} \my(\lambda_n, H(\ccdot, u_n))(z)\\
&\leq \limsup_{n \to \infty} \my(\lambda_n,H(\ccdot,u_n))(z) \leq \my(\lambda,H(\ccdot,u))(z).
\end{eq*}
Thus, we deduce the pointwise convergence and by \cite[Theorem 3.26]{bib:AttouchEpi} the asserted Mosco-convergence.
\end{proof}
In the light of the obstacle problem \refer{eq:obstacleproblem} and \cite{bib:AlphonseHintRautenbergDirectDiff} an important class of \msa{applications for the previous few lemmas}{} are indicator functionals of a set-valued mapping $K_{L^2}(u) := \{z \in L^2(\Omega): z \geq \Phi(u)\}$, 
where $U \overset{c}{\hookrightarrow} L^2(\Omega) \msa{=: Z}{}$ embeds compactly. \msa{In this setting, take}{} $H(z,u):= I_{K_{L^2}(u)}(z)$. \msa{Then, $H$}{} fulfills the conditions of \refer{lem:my:mosco} and \refer{thm:myreg:gammaconv} with $B = i_{L^2(\Omega)}$:\\
Take sequences $z_n \rightharpoonup z$ in $L^2(\Omega)$, $v_n \rightharpoonup v$ and $u_n \rightharpoonup u$ in $U$ and assume $\Phi:U \to U$ to be \msa{\emph{weakly continuous}}{}, i.e. that for $u_n \rightharpoonup u$ holds $\Phi(u_n) \rightharpoonup \Phi(u)$ in $U$. Then, we obtain by the weak closedness of the set of all non-negative $L^2(\Omega)$-functions, that \msa{$z_n \in K_{L^2}(u_n)$ and $v_n \in K_{L^2}(u_n) \cap U$ imply}{}
\begin{eq*}
z_n - \Phi(u_n) \rightharpoonup z - \Phi(u) \geq 0 \text{ and } v_n - \Phi(u_n) \rightharpoonup v - \Phi(u) \geq 0.
\end{eq*}
On the other hand, for given $u_n \rightharpoonup u$ and $z \in K_{L^2}(u)$ and $v \in K_{L^2}(u) \cap U$ we construct the sequences $z_n = z - \Phi(u) + \Phi(u_n)$ and $v_n = v - \Phi(u) + \Phi(u_n)$. By construction holds $z_n \in K_{L^2}(u_n)$ and $v_n \in K_{L^2}(u_n) \cap U$ and using the weak continuity of $\Phi$ we deduce \msa{$v_n \rightharpoonup v$ in $U$ and $z_n \to z$ in $Z$ by the compact embedding of $U$ into $Z$}{}. This guarantees the conditions on the lower and upper limits in \refer{thm:myreg:gammaconv} and \refer{lem:my:conv}.\\
After establishing the $\Gamma$-convergence for the regularized part of the functional, we are now ready to prove our key convergence result \refer{thm:myreg:gammaconv}.

\begin{proof}[Proof of \refer{thm:myreg:gammaconv}]
Let an arbitrary sequence $u_n \rightharpoonup u$ be given. By \refer{lem:my:mosco} we obtain the Mosco-convergence
\begin{eq*}
\my(\gamma_n,H(\ccdot,u_n)) \moscoconv H(\ccdot,u).
\end{eq*}
Hence, we obtain by the weak lower semi-continuity and for $Bu_n \rightharpoonup Bu$ that
\begin{eq*}
\calE(u,u) &= F(u,u) + G(Au,u) + H(Bu,u) \leq \liminf_{n \to \infty} F(u_n,u_n)\\
&+ \liminf_{n \to \infty} G(Au_n,u_n) + \liminf_{n \to \infty} \my(\gamma_n, H(\ccdot,u_n))(Bu_n)\\
&\leq \liminf_{n \to \infty} \calE_n(u_n,u_n).
\end{eq*}
Take an arbitrary $v \in U_\ad$. To construct a recovery sequence, we take $v_n \rightharpoonup v$ to be the recovery sequence such that
\begin{eq*}
\limsup_{n \to \infty} \my(\gamma_n, H(\ccdot, u_n))(Bv_n) \leq H(Bv,u).
\end{eq*}
By the continuity assumptions on $F$ and $G$ we obtain
\begin{eq*}
&\limsup_{n \to \infty} \left( F(v_n,u_n) \right.\left.+ G(A v_n,u_n) + \my(\gamma_n,H(\ccdot,u_n))(Bv_n) \right)\\
&\leq \limsup_{n \to \infty} F(v_n,u_n) + \limsup_{n \to \infty} G(A v_n,u_n)\\
&+ \limsup_{n \to \infty} \my(\gamma_n, H(\ccdot,u_n))(Bv_n) \leq F(v,u) + G(Av,u) + H(Bv,u).
\end{eq*}
\msa{This proves in combination the asserted weak $\Gamma$-convergence.}{}
\end{proof}
In a similar fashion, one can prove a slightly stronger result guaranteeing Mosco-convergence in the case of $H$ being independent of the feedback component.
\begin{thm}\label{thm:myreg:mosco}
Let proper, convex, lsc. functionals $F:U_\ad \times U_\ad \to \bR, G: Y \times U_\ad \to \bR$ be given as in \refer{thm:myreg:gammaconv} and let $H:Z \to \bR$ be given as in \refer{lem:my:gamma}. Define the functionals
\begin{eq*}
\calE(v,u) &:= F(v,u) + G(Av,u) + H(Bv) \text{ and}\\
\calE_n(v,u) &:= F(v,u) + G(Av,u) + \my(\gamma_n, H)(Bv).
\end{eq*}
Then, the convergence $\calE_n \moscoconv \calE$ holds.
\end{thm}
\begin{proof}
The proof uses \refer{lem:my:gamma} \msa{and}{} is analogous to the one in \refer{thm:myreg:gammaconv} and is therefore omitted.
\end{proof}
So far, we only treated the $\Gamma$-convergence of Moreau--Yosida regularized problems in the infinite-dimensional function space setting. As thus, no discretization has been imposed. This will be addressed in the next subsection.
\subsection{$\Gamma$-Convergence for Regularized and Dis\-cre\-tized Equilibrium Problems}
There are several successful discretization methods discussed in the numerical literature, among them the most prominent ones in the context of partial differential equations are \emph{finite difference methods (FDM)} and \emph{finite element methods (FEM)}. We devote our attention to the latter. The general idea for discretizing optimization problems is to restrict the minimization to a finite dimensional subspace. For equilibrium problems this procedure works analogously and an abstract result is given in the following Lemma.
\begin{lem}\label{lem:fem:mosco}
Let a functional \msa{$\calE: U_\ad \times U_\ad \to \bR$}{} be given, such that for all sequences $u_\ell \rightharpoonup u$ holds
\begin{eq*}
\calE(u,u) \leq \liminf_{\ell \to \infty} \calE(u_\ell,u_\ell) 
\end{eq*}
and for all sequences $u_\ell \rightharpoonup u$ and $v_\ell \to v$ holds
\begin{eq*}
\limsup_{\ell \to \infty} \calE(v_\ell,u_\ell) \leq \calE(v,u).
\end{eq*}
Let a sequence of subspaces $(U_\ell)_{\ell \in \N}$ be given with $U_\ell \subseteq U_{\ell+1} \subseteq U$ and $U_\ad \cap U_\ell \neq \emptyset$ as well as $\dom{\calE(\ccdot,u_\ell)} \neq \emptyset$ for all $u_\ell \in U_\ad \cap U_\ell$ and $\ell \in \N$, such that
\begin{eq*}
U = \cl{U_\ad \cap \bigcup_{\ell \in \N} U_\ell}
\end{eq*}
and define the functionals $\calE_\ell: U_\ad \times U_\ad \to \bR$ by 
\msa{\begin{eq*}
\calE_\ell(v,u) := \calE(v,u) + I_{U_\ell}(v).
\end{eq*}}%
Then, the Mosco-convergence $\calE_\ell \moscoconv \calE$ holds true.
\end{lem}
\begin{proof}
Let a sequence $u_\ell \msa{\rightharpoonup}{} u$ be given. On the one hand, we obtain by assumption
\begin{eq*}
\calE(u,u) \leq \liminf_{\ell \to \infty} \calE(u_\ell,u_\ell) \leq \liminf_{\ell \to \infty} \calE_\ell(u_\ell,u_\ell),
\end{eq*}
since $0 \leq I_{U_\ell}(u_\ell)$. On the other hand, take an arbitrary $v \in U_\ad$. Then, by the assumed density there exists a sequence $v_\ell \to v$ with $v_\ell \in U_\ad \cap U_\ell$ and we obtain
\begin{eq*}
\limsup_{\ell \to \infty} \calE_\ell(v_\ell,u_\ell) = \limsup_{\ell \to \infty} \calE(v_\ell,u_\ell) \leq \calE(v,u),
\end{eq*}
which proves the Mosco-convergence.
\end{proof}
The monotonicity condition imposed on the subspaces associates to a sequence of triangulations constructed via successive refinements of the mesh. The density condition on the other hand translates as the longest edge of the triangulation going to zero.\\
The methodology is to first apply the Moreau--Yosida regularization to the original problem and afterwards, impose the discretization. In the previous subsection, we discussed the behaviour with respect to the regularization only. Next, we apply \refer{lem:fem:mosco} to obtain a result with respect to the  discretization only.
\begin{thm}
Let the functionals $F,G,H$ be given as in \refer{thm:myreg:gammaconv} (respectively \refer{thm:myreg:mosco}) along with a sequence of subspaces $(U_\ell)_{\ell \in \N}$ with $U_\ell \subseteq U_{\ell+1} \subseteq U$ with $U_\ad \cap U_\ell \neq \emptyset$ as well as $\dom{F(\ccdot,u_\ell) + G(A \ccdot,u_\ell) } \neq \emptyset$ for all $u_\ell \in U_\ad \cap U_\ell$ and
\begin{eq*}
U = \cl{U_\ad \cap \bigcup_{\ell = 1}^\infty U_\ell}.
\end{eq*}
Define \msa{for fixed $\gamma > 0$}{} the following functionals
\begin{eq*}
\calE(v,u) &:= F(v,u) + G(Av,u) + \my(\gamma, H(\ccdot, u))(Bv)\text{ and}\\
\calE_\ell(v,u) &:= \calE(v,u) + I_{U_\ell}(v).
\end{eq*}
Then, the convergence $\calE_\ell \moscoconv \calE$ holds true.
\end{thm}
\begin{proof}
In the setting of \refer{thm:myreg:gammaconv} we check the conditions of \refer{lem:fem:mosco}. \msa{For that, we}{} take a sequence $u_n \rightharpoonup u$ and \msa{obtain}{} by using \refer{lem:my:conv} the estimate
\begin{eq*}
F(u,u) &+ G(Au,u) + \my(\gamma,H(\ccdot, u))(Bu) \leq \liminf_{\ell \to \infty} F(u_\ell,u_\ell)\\
&+ \liminf_{\ell \to \infty} G(Au_\ell,u_\ell) + \liminf_{\ell \to \infty} \my(\gamma, H(\ccdot,u_\ell))(B u_\ell)
\end{eq*}
and for every sequence $v_\ell \to v$ the result
\begin{eq*}
&\limsup_{\ell \to \infty} \left(F(v_\ell,u_\ell) + G(Av_\ell,u_\ell) + \my(\gamma,H(\ccdot,u_\ell))(B u_\ell) \right)\\
&\leq \limsup_{\ell \to \infty} F(v_\ell,u_\ell) + \limsup_{\ell \to \infty} G(Au_\ell,u_\ell) + \limsup_{\ell \to \infty} \my(\gamma,h(\ccdot,u_\ell))(Bu_\ell)\\
&\leq F(v,u) + G(Av,u) + \my(\gamma,H(\ccdot,u))(Bv).
\end{eq*}
\end{proof}
As one is able to see, the upper semi-continuity of the equilibrium functional was relevant for the proof. Hence, the application of the Moreau--Yosida regularization in first place was vital for the application of \refer{lem:fem:mosco}.\\
So far, we have discussed the $\Gamma$-convergence of \msa{the}{} Moreau--Yosida regularization as well as of \msa{its}{} discretizations. What remains is the discussion of the behaviour when both is applied simultanously. To obtain a practically relevant convergence result, we propose the following Theorem.
\begin{thm}[Diagonal theorem]\label{thm:diag}
Let families of functionals $(\calE_n^\ell)_{\ell,n\in\N}$, $(\calE_n)_{n \in \N}, \calE : U \times U \to \bR$ and sequences $(\msa{u^\ell_n}{})_{\ell,n\in\N}, (u_n)_{n \in \N} \subset U$ be given along with $u \in U$ as well as \msa{sequences}{} $(\eps_n^\ell)_{\ell,n\in \N},(\eps_n)_{n \in \N} \subset (0,\infty)$ with $\eps^\ell_n \to 0$ and $\eps_n \to 0$ such that 
\begin{eq*}
\calE_n^\ell(u_n^\ell,u_n^\ell) \leq \inf \calE_n^\ell(\ccdot,u_n^\ell) + \eps_n^\ell \text{ and } \calE_n(u_n,u_n) \leq \inf \calE_n(\ccdot,u_n) + \eps_n
\end{eq*}
\msa{hold true}{}. Let the following convergence assumptions
\begin{eq*}
\calE_n^\ell \moscoconv \calE_n \text{ and } u_n^\ell \rightharpoonup u_n \text{ as } \ell \to \infty 
\end{eq*}
and
\begin{eq*}
\calE_n \moscoconv \calE \text{ and } u_n \rightharpoonup u \text{ as } n \to \infty
\end{eq*}
hold. Then, there exists a subsequence $\ell_n \to \infty$ such that
\begin{eq*}
\inf \calE(\,\cdot\,, u) = \inf \calE_n^{\ell_n}(\,\cdot\,,u_n^{\ell_n})
\end{eq*}
Moreover, if $\left(\calE^\ell_n\right)_{\ell, n \in \N}, (\calE_n)_{n \in \N}, \calE$ are strongly convex with the same constant, then the sequence $\ell_n$ can be chosen, such that $u_n^{\ell_n} \to u$ converges.
\end{thm}
\begin{proof}
\msa{First}, we observe \msa{by}{} the choice of $u_n$ and the \msa{assumed}{} Mosco-convergence \msa{using}{} a recovery sequence $v_n \to v$ \msa{with}{} $\limsup_{n \to \infty} \calE_n(v_n,u_n) \leq \calE(v,u)$ \msa{that}
\begin{eq*}
\calE(u,u) &\leq \liminf_{n \to \infty} \calE_n(u_n,u_n) \leq \liminf_{n \to \infty} \left( \inf \calE_n(\ccdot,u_n) + \eps_n \right)\\
&= \liminf_{n \to \infty} \msa{\left(} \inf \calE_n(\ccdot,u_n) \msa{\right)} \leq \limsup_{n \to \infty} \msa{\left(} \inf \calE_n(\ccdot,u_n) \msa{\right)}\\
&\leq \limsup_{n \to \infty} \calE_n(v_n,u_n) \leq \calE(v,u).
\end{eq*}
Since $v \in U$ was chosen arbitrarily we observe
\begin{eq*}
\inf \calE(\ccdot,u) \leq \liminf_{n \to \infty} \msa{\left(\inf}{} \calE_n(\ccdot,u_n) \msa{\right)} \leq \limsup_{n \to \infty} \msa{\left( \inf}{} \calE_n(\ccdot, u_n) \msa{\right)} \leq \inf \calE(\ccdot,u)
\end{eq*}
and hence 
\begin{eq}[\label{eq:tmp:diagthm}]
\lim_{n \to \infty} \left( \inf\calE_n(\ccdot, u_n) \right) = \inf \calE(\ccdot,u) = \calE(u,u).
\end{eq}
With \msa{an analogous}{} estimate \msa{carried out for the Mosco-convergence $\calE_n^\ell \moscoconv \calE_n$ as $\ell \to \infty$}{} we deduce
\begin{eq*}
\inf \calE_n(\ccdot,u_n) &\leq \liminf_{n \to \infty} \left(\msa{\inf}{} \calE_n^\ell(\ccdot,u^\ell_n) +\eps_n^\ell \right)\text{ and }\\
\limsup_{n \to \infty} \left(\inf \calE_n^\ell(\ccdot,u^\ell_n) \right) &\leq \inf \calE_n(\ccdot,u_n).
\end{eq*}
By choosing for every $n \in \msa{\N}$ a natural number $\ell_n$ such that $\eps_n^\ell \leq 2\eps_n$ and
\begin{eq}[\label{eq:tmp2:diagthm}]
\inf \calE_n(\ccdot,u_n) &\leq \inf \calE_n^\ell(\ccdot,u_n^\ell) + \eps_n^\ell + \eps_n \text{ as well as}\\
\inf \calE_n^\ell(\ccdot,u_n^\ell) &\leq \inf \calE_n(\ccdot,u_n) + \eps_n
\end{eq}
for all $\ell \geq \ell_n$\msa{,}{} we obtain by the repeated use of \refer{eq:tmp:diagthm} \msa{and \refer{eq:tmp2:diagthm}}{} that
\begin{eq*}
\inf \calE(\ccdot,u) &= \calE(u,u) = \lim_{n \to \infty} \left(\inf \calE_n(\ccdot,u_n)\right)\\
&\leq \liminf_{n \to \infty} \left(\inf \calE_n^{\ell_n}(\ccdot,u_n^{\ell_n}) + 3\eps_n\right)\\
&= \liminf_{n \to \infty} \left( \inf \calE_n^{\ell_n}(\ccdot,u_n^{\ell_n}) \right) \leq \limsup_{n \to \infty} \left( \inf \calE_n^{\ell_n}(\ccdot,u_n^{\ell_n}) \right)\\
&\leq \limsup_{n \to \infty} \left(\inf \calE_n(\ccdot,u_n) + \eps_n \right) \leq \inf \calE(\ccdot,u).
\end{eq*}
Thus, we deduce the first assertion.\\
For the second assertion we take for every $n \in \N$ a recovery sequence $v_n^\ell \to u_n$ with $\limsup_{\ell \to \infty} \calE_n^\ell(v_n^\ell,u_n^\ell) \leq \calE_n(u_n,u_n)$ as $\ell \to \infty$ as well as a recovery sequence $v_n \to u$ such that $\limsup_{n \to \infty} \calE_n(v_n,u_n)\leq \calE(u,u)$. Then, we deduce \msa{the}{} following estimate
\begin{eq*}
\|u-u_n^\ell\| &\leq \|u - v_n\| + \|v_n - \calB_{\msa{n}}(u_n)\| + \|\calB_n(u_n) - u_n\|\\
&+ \|u_n - v_n^\ell\| + \|v_n^{\msa{\ell}} - \calB_n^\ell(u_n^\ell)\| + \|\calB_n^\ell(u_n^\ell) - u_n^\ell\|.
\end{eq*}
Let $\ell_n'$ denote the previously constructed sequence and choose $\ell_n''$ such that
\begin{eq*}
\calE_n^\ell(v_n^\ell,u_n^\ell) \leq \calE_n(u_n,u_n) + \eps_n \text{ and } \|v_n^\ell - u_n\| \leq \eps_n
\end{eq*}
for all $\ell \geq \ell_n''$. Then, we deduce for $\ell\geq \max(\ell_n',\ell_n'')$
\begin{eq*}
\frac{\alpha}{2}\|v_n^\ell &- \calB_n^\ell(u_n^\ell)\|^2 \leq \calE_n^\ell(v_n^\ell,u_n^\ell) - \inf \calE_n^\ell(\ccdot,u_n^\ell)\\
&\leq \calE_n^\ell(v_n^\ell,u_n^\ell) - \inf \calE_n(\ccdot,u_n) + \eps_n + \eps_n^\ell\\
&\leq \calE_n(u_n,u_n) -\inf \calE_n(\ccdot,u_n) + 4\eps_n \leq 5\eps_n
\end{eq*}
as well as
\begin{eq*}
\frac{\alpha}{2}\|\calB_n^\ell(u_n^\ell) - u_n^\ell\|^2 \leq \calE_n^\ell(u_n^\ell,u_n^\ell) - \inf \calE_n^\ell(\ccdot,u_n^\ell) \leq \eps_n^\ell \leq 2\eps_n
\end{eq*}
by the assumption on $u_n^{\msa{\ell}}$. Analogously\msa{,}{} we get
\begin{eq*}
\frac{\alpha}{2}\|u_n - \calB(u_n)\|^2 \leq \calE_n(u_n,u_n) - \inf \calE_n(\ccdot,u_n) \leq \eps_n
\end{eq*}
as well as 
\begin{eq*}
\frac{\alpha}{2}\|v_n-\calB_n(u_n)\|^2 \leq \calE_n(v_n,u_n) - \inf \calE_n(\ccdot,u_n).
\end{eq*}
Hence, we deduce in total with $\ell_n := \max(\ell_n',\ell_n'')$ 
\begin{eq*}
0 \leq \limsup_{n \to \infty}\|u&-u_n^{\ell_n}\| \leq \limsup_{n \to \infty}\left( \|u-v_n\| + \calE_n(v_n,u_n)\right.\\
&\left.- \inf\calE_n(\ccdot,u_n) + 8\eps_n \right) \leq \calE(u,u) - \inf\calE(\ccdot,u) = 0
\end{eq*}
and hence the convergence $u_n^{\ell_n} \to u$.
\end{proof}
The Diagonal theorem tells us, that \msa{given}{} a hierarchical sequence of $\Gamma$-convergent functionals \msa{along with a sequence of approximate equilibria}{}, we can select a subsequence \msa{of the latter}{} converging to a solution of the \msa{limiting}{} problem. The convergence result \refer{prop:equi:gamma:mini_convergence} can be interpreted as a corollary of \refer{thm:diag}.\\
With these tools at hand we close our discussion of a priori results on convergence and turn our attention to adaptive refinement techniques. 
\section{A Posteriori Theory}
In order to obtain a reliable error estimator we generalize the techniques in \refer{bib:BartelsMilicevicPrimalDualGap}. Therein, techniques for convex minimization problems have been developed and applied to a selection of problems structurally related \msa{to the}{} problems under investigation in this paper. There, the error is the norm of the difference between the approximation and the unique solution has been estimated. In general, the existence of equilibria can be a challenging problem. Moreover, even in seemingly simple cases\msa{,}{} uniqueness does not need to hold \msa{(cf. \cite[Example 1]{bib:KanzowFacchinei})}{}. As thus, using the distance to the solution set is a challenging object to estimate. Instead, we seek to estimate a meaningful residuum. For this sake, we \msa{introduce the}{} \emph{best response operator}.
%
%
\begin{defn}[\msa{Best}{} response operator]
Let a functional $\calE: U_\ad \times U_\ad \to \bR$ with $\dom{\calE(\cdot,u)} \neq \emptyset$ be given. The \emph{best response operator} $\calB: U_\ad \rra U_\ad$ is defined \msa{by}{}
\begin{eq*}
\calB(u) := \argmin_{v\in U_\ad} \calE(v,u).
\end{eq*}
\end{defn}%
\msa{Evidently,}{} a point is an equilibrium of $\calE : U_\ad \times U_\ad \to \bR$ with best response mapping $\calB: U_\ad \rra U_\ad$\msa{,}{} if and only if 
\begin{eq*}
u \in \calB(u)
\end{eq*}
holds.
This observation \msa{is}{} used to formulate a residuum. In the light of \cite[Proposition 3.1]{bib:BartelsMilicevicPrimalDualGap} this motivates the following Theorem.

\begin{thm}\label{thm:estimate}
Let a non-empty, convex, closed subset \msa{$U_\ad$}{} of some Banach space $U$ \msa{and}{} a functional $\calE:U_\ad \times U_\ad \to \bR$ \msa{with $\calE(\cdot,u)$ being \lsc{} and $\dom{\calE(\cdot,u)} \neq \emptyset$ for all $u \in U_\ad$,}{} and assume that there exists $\alpha > 0$ such that $\calE(\ccdot,u):U_\ad \to \bR$ is 
\emph{$\alpha$-strongly}-convex for all $u \in U_\ad$.\\
Then, the \msa{best response}{} operator $\calB$ is singleton and for all $u \in U_\ad$ holds
\begin{eq*}
\frac{\alpha}{2}\| v - \calB(u) \|_{\msa{U}}^2 \leq \calE(v,u) - \calE(\calB(u),u).
\end{eq*}
\end{thm}
\begin{proof}
Let arbitrary $v, u \in U_\ad$ be given. Since \msa{by assumption}{} $\calE(\ccdot, u) : U_\ad \rightarrow \bR$ is \msa{proper, lsc. and }{} strongly convex, \msa{it is in particular strictly convex and thus}{} the best response mapping is singleton. \msa{By testing the control component with $tv + (1-t)\calB(u)$ for an arbitrary $t \in (0,1)$}{} we deduce from the strong convexity
\begin{eq*}
&\scalebox{0.948}{$\dfrac{\alpha}{2} t(1-t)\|v-\calB(u)\|_{\msa{U}}^2 \leq t\calE(v,u) + (1-t)\calE(\calB(u),u)- \calE(tv + (1-t)\calB(u),u)$}\\
&\scalebox{0.948}{$\hphantom{\dfrac{\alpha}{2} t(1}= t(\calE(v,u) - \calE(\calB(u),u)) + \calE(\calB(u),u) - \calE(tv + (1-t)\calB(u),u)$}\\
&\scalebox{0.948}{$\hphantom{\dfrac{\alpha}{2} t(1}\leq t(\calE(v,u) - \calE(\calB(u),u))$}.
\end{eq*}
Dividing by $t > 0$ and passing $t \searrow 0$ yields the assertion.
\end{proof}
In practice\msa{,}{} $u$ plays the role of an approximated solution and as such we are in practice interested in the case $v = u$.\\
It is worth noting, that the requirement of $\calE$ being strongly convex in the control component is less restrictive for equilibria, than it is for optimization problems\msa{:}{} So see this, consider a functional $\calE: U_\ad \times U_\ad \to \bR$ being \msa{proper, \lsc{} and}{} convex in its control component and define a modified functional using a parameter $\alpha > 0$ \msa{by}{}
\begin{eq*}
\calE_\alpha(v,u) := \calE(v,u) + \frac{\alpha}{2}\|v - u\|_{\msa{U}}^2,
\end{eq*}
see also \cite{bib:KanzowHeusingerOptReformulation, bib:JSPangGuerkan}. Then $\calE_\alpha: U_\ad \times U_\ad$ is \msa{$\alpha$-}strongly convex in the control component.  In fact, both formulations have the same set of equilibria as it can be seen as follows:\\
Let first $u \in U_\ad$ be an equilibrium of $\calE$. Then\msa{,}{} we obtain for all $ v\in U_\ad$ the inequality
\begin{eq*}
\calE_\alpha(u,u) = \calE(u,u) \leq \calE(v,u) \leq \calE_\alpha(v,u)
\end{eq*}
and hence $u$ being an equilibrium of $\calE_\alpha$.\\
\msa{For the other direction,}{} let now an equilibrium $u \in U_\ad$ of the functional $\calE_\alpha$ be given. Take arbitrary $v \in U_\ad$ and $t \in (0,1)$ and we obtain
\begin{eq*}
\calE(u,u) &= \calE_\alpha(u,u) \leq \calE_\alpha(tv + (1-t)u,u)\\
&= \calE(t v + (1-t) u, u) + \frac{\alpha}{2}\|tv + (1-t)u - u\|_H^2\\
&\leq t \calE(v,u) + (1-t)\calE(u,u) + \frac{\alpha t^2}{2}\|v - u\|_H^2\msa{.}
\end{eq*}
Subtracting $(1-t)\calE(u,u)$ on both sides and dividing by $t$ yields
\begin{eq*}
\calE(u,u) \leq \calE(v,u) + \frac{\alpha t}{2}\|v-u\|_H^2.
\end{eq*}
Passing $t \searrow 0$ proves $u$ being an equilibrium of $\calE$.\\
Since the \msa{underlying}{} minimization is \msa{a}{} difficult \msa{task}, we seek to dualize the optimization problem $\min_{v \in U_\ad} \calE(v,u)$ for given $u \in U_\ad$. Therefore, we introduce the notion of a dual functional in the context of equilibrium problems. 
\begin{defn}\label{defn:dual}
Let a functional $\calE:U_\ad \times U_\ad \rightarrow \bR$ be given and let $\calP$ denote a set. A functional $\calD: \calP \times U_\ad \rightarrow \bR$ is called a \emph{dual functional} of $\calE$, if for all $v,u \in U_\ad$ and all $p \in \calP$ holds
\begin{eq*}
0 \leq \calE(v,u) + \calD(p,u).
\end{eq*}
Moreover, $\calD$ is called a \emph{strong dual functional}, if
\begin{eq*}
\min_{v \in U_\ad} \calE(v,u) + \inf_{p \in \calP} \calD(p,u) = 0 \text{ for all } u \in U_\ad
\end{eq*}
holds.
\end{defn}
\msa{By definition}, every strong dual functional is a dual functional. For every \msa{strong dual functional}{} holds \msa{under the assumptions of}{} \refer{thm:estimate} for all $u \in U_\ad$ and $p \in \calP$ the estimate
\begin{eq}\label{eq:aposteriori}
\frac{\alpha}{2}\|u - \calB(u)\|_U^2 \leq \eta^2(p,u) := \calE(u,u) + \calD(p,u).
\end{eq}
For given $v,u \in U_\ad$, the sharpness of the estimate depends solely on the choice of $p \in \calP$, which can be chosen arbitrarily.\\
It is worth noting, that in the light of the Fenchel duality theorem, given \msa{in}{} \refer{thm:fencheldual} in the appendix, the roles of the equilibrium functional and the dual functional are swapped: The minimum is supposed to be achieved with the variable $u$, but with respect to the dual variable $p$ only the infimum is needed. Hence, for the derivation of a strong dual functional the formulation of a \msa{\emph{pre-dual problem}}{} might be of interest. We will briefly return to that point towards the end of this section.\\
Next, we derive for problems \msa{of the type presented in}{} \refer{eq:intro} \msa{an associated}{} dual functional. \msa{Since the dualization takes place with respect to the control component, the following result is formulated for optimization problems.}{}
\begin{lem}\label{lem:duality:my}
Let Banach spaces $U,Y$ and a Hilbert space $Z$\msa{, with its dual space $Z^*$ identified with $Z$,}{} \msa{as well as convex}{}, proper, \lsc{} functionals $f: U \to \bR$, $g:Y \to \bR$ and $h:Z \to \bR$ \msa{be given.}{} \msa{Moreover, let bounded}, linear operators $A \in \calL(U,Y)$ and $B \in \calL(U,Z)$ and a positive real number $\gamma > 0$ \msa{be given}{}. Consider the following minimization problem
\begin{eq}\label{eq:duality:my:primal}
\inf_{u \in \msa{U}} \left( f(u) + g(Au) + \my(\gamma,h)(Bu)\right).
\end{eq}
If the following \emph{constraint qualification}
\begin{eq}\label{eq:duality:my:cq}
0 \in \core{ \dom g - A \dom f}
\end{eq}
is fulfilled, then the Fenchel dual problem \msa{of \refer{eq:duality:my:primal}}{} reads as
\begin{eq}\label{eq:duality:my:dual}
\inf_{y^* \in Y^*, z^* \in \msa{Z}} \left( f^*(-A^*y^* - B^*z^*) + g^*(y^*) + h^*(z^*) + \frac{1}{2\gamma}\|z^*\|^2_{\msa{Z}}\right).
\end{eq}
\end{lem}
\begin{proof}
We check the conditions of \refer{thm:fencheldual}. Rewriting \refer{eq:duality:my:primal} using \refer{defn:my} leads to 
\begin{eq}\label{eq:duality:my:primal:reform}
\inf_{u \in U, \zeta \in Z} \left( f(u) + g(Au) + h(\zeta) + \frac{\gamma}{2}\|\zeta - Bu\|_Z^2 \right).
\end{eq}
Introducing the functionals $\widetilde f: U \times Z \to \bR$ defined by $\widetilde f (u,\zeta) := f(u) + h(\zeta)$ and $\widetilde g: Y \times Z \times Z \to \bR$ defined by $\widetilde g(y,z) := g(y) + \frac{\gamma}{2}\|z\|_Z^2$ \msa{as well as}{} the linear continuous operator $\widetilde A \in \calL(U \times Z, Y \times Z)$ defined by $\widetilde A(u,\zeta) = (Au,\zeta - Bu)$ \msa{we}{} rewrite \refer{eq:duality:my:primal:reform} as 
\begin{eq*}
\inf_{(u,\zeta) \in U \times Z} \widetilde f(u,\zeta) + \widetilde g\left(\widetilde A(u,\zeta)\right).
\end{eq*}
In the light of the Fenchel duality theorem in \refer{thm:fencheldual}, given in the appendix, the constraint qualification \refer{eq:fencheldual:cq} reads as
\begin{eq*}
0 &\in \core{\dom{\widetilde g} - \widetilde A \dom{\widetilde f}}\\
&= \core{\dom{g} \times Z - (A\dom f) \times (\dom h - B \dom f) },
\end{eq*}
which is equivalent to 
\begin{eq*}
0 \in \core{\dom{g} - A \dom f}
\end{eq*}
as proposed in the \msa{assumption}. The calculation of the con\-ju\-gate functionals \msa{yields}{} 
\begin{eq*}
\widetilde f^*(u^*,\zeta^*) = f^*(u^*) + h^*(\zeta^*),
\end{eq*}
as well as 
\begin{eq*}
\widetilde g^*(y^*, z^*) = g^*(y^*) + \frac{1}{2\gamma}\|z^*\|_Z^2\msa{.}
\end{eq*}
The dual operator $\widetilde A^*$ reads as $\widetilde A^* (y^*, z^*) = (A^* y^* - B^* z^*, z^*)$. Hence\msa{,}{} the dual problem can be formulated as
\begin{eq*}
\min_{y^*\in Y^*, z^* \in Z} \left( \widetilde f^*(- \widetilde A^* (y^*, z^*)) + \widetilde g^*(y^*, z^*) \right)
\end{eq*}
\msa{or equivalently}{}
\begin{eq*}
\min_{y^*\in Y^*, z^* \in Z} \left( f^*(-A^* y^* - B^* z^*) + h^*(z^*) + g^*(y^*) + \frac{1}{2\gamma}\|z^*\|_Z^2\right),
\end{eq*}
which yields the assertion.
\end{proof}
For \msa{the constraint qualification}{} the domains of the functionals play a role. As the Moreau--Yosida regularization has \msa{as}{} its domain the whole space, the domain of $h$ itself does not influence the constraint qualification. Thus, one might be able to dualize the regularized problem, but not the original one.\\
Returning to the obstacle problem \refer{eq:obstacleproblem} we can formulate the dual functional as follows: As derived in the discussion to \refer{defn:my} the regularized equilibrium functional reads as
\begin{eq*}
\calE_\gamma(y) := \frac{1}{2}\|\nabla y\|_{L^2(\Omega;\R^d)}^2 - (f, i_{L^2}y) + \frac{\gamma}{2}\int_\Omega (\lowobs - y)^{2+}\dx.
\end{eq*}
To apply \refer{lem:duality:my} we decompose the problems as follows\msa{:}{}\\
\msa{Set for the spaces}{} $U = H_0^1(\Omega)$, $Y = L^2(\Omega;\R^d) \times L^2(\Omega)$ and $Z = L^2(\Omega)$. \msa{For the functionals}{}  $F: H_0^1(\Omega) \to \bR$ with $F(y) = 0$ for all $y \in H_0^1(\Omega)$, and $G: Y \to \bR$ with $G(p,z) := \frac{1}{2}\|p\|_{L^2(\Omega;\R^d)}^2 - (f,z)_{L^2(\Omega)}$ together with $H:Z \to \bR$, $H(z):= I_{K_{L^2}}(z)$ \msa{are chosen.}{} \msa{For the}{} operators $A \in \calL(U, Y)$ defined by $Ay := (\nabla y, i_{L^2(\Omega)} y)$ and \msa{$B \in \calL(U,Z)$,}{} $B = i_{L^2(\Omega)}$ \msa{are taken}{}.\\
By \refer{lem:duality:my} together with the results of Calculation \ref{calc:conjugate} \msa{in the appendix,}{} we deduce after some slight simplifications the strong dual functional
\begin{eq*}
\calD_\gamma(p,z) &:= \frac{1}{2}\|p\|_{L^2(\Omega; \R^d)}^2 + (z,\lowobs)_{\msa{L^2(\Omega)}} + \frac{1}{2\gamma}\|z\|_{L^2(\Omega)}^2\\
&\ + I_{\{0\}}(-\ddiv p + z - f) + I_{L^2_-(\Omega)}(z)\msa{,}
\end{eq*}
\msa{with $L^2_-(\Omega) := \{z \in L^2(\Omega) : z \leq 0 \text{ a.e. on }\Omega \}$.}{}\\
In the light of \msa{the previous convergence discussion for regularized problems in \refer{sec:regularizedgamma} one}{} might ask as well for the behavior of the dual functionals \msa{with respect to the inherited regularization parameter $\gamma$}{}. In fact, for $\gamma \to \infty$ one expects to find \msa{as}{} limiting functional 
\begin{eq*}
\calD(p,z) := \frac{1}{2}\|p\|_{L^2(\Omega; \R^d)}^2 + (z,\lowobs)_{\msa{L^2(\Omega)}} + I_{\{0\}}(-\ddiv p + z - f) + I_{L^2_-(\Omega)}(z).
\end{eq*}
\msa{Thus, one}{} might ask for the relation between the original functional $\calE$ and $\calD$. Returning to \msa{\refer{defn:dual}}, we can integrate the limiting dual functional $\calD$ as the \emph{pre-dual} of $\calE$. To show this, we use the Fenchel duality theorem (see \refer{thm:fencheldual}) 
with $\msa{U}{} := L^2(\Omega;\R^d) \times L^2(\Omega)$ and $Y := H^{-1}(\Omega) = (H_0^1(\Omega))^*$ along with $f : \msa{U}{} \to \bR$ with $f(p,z) := \frac{1}{2}\|p\|_{L^2(\Omega;\R^d)}^2 + (z,\lowobs)_{L^2(\Omega)} + I_{L^2_-(\Omega)}(z)$ as well as $g:Y \to \bR$ with $g(y):= I_{\{f\}}(y)$ and $A \in \calL(X,Y)$ with $A(p,z) := -\ddiv p + i_{H^{-1}(\Omega)}z$ to it and obtain as constraint qualification 
\begin{eq*}
0 &\in \core{\dom g - A \dom f}\\
&= \core{ \{f\} + \ddiv L^2(\Omega;\R^d) - i_{H^{-1}(\Omega)} L^2(\Omega)}.
\end{eq*}
\msa{For its proof}, it is sufficient to \msa{show}{} that the set inside the core is indeed the space $H^{-1}(\Omega)$. Therefore\msa{, for given $\xi \in H^{-1}(\Omega)$}{} solve the partial differential equation 
\begin{eq*}
-\Delta y = f - \xi \text{ in } \Omega,\ y = 0 \text{ on } \partial \Omega.
\end{eq*}
\msa{Clearly, it has a unique solution $y \in H_0^1(\Omega)$}{} and \msa{setting}{} $p = \nabla y \in L^2(\Omega;\R^d)$ as well as $z = 0 \in L^2(\Omega)$ \msa{yields}{} $\xi = f + \ddiv p$ in $H^{-1}(\Omega)$, which proves the assertion.\\
Under the use of Calculation \ref{calc:conjugate} and the Fenchel biconjugation theorem (see e.g. \cite[Theorem 4.2.1]{bib:JBorweinLewisConvexAnalysis}) we recover after slight simplifications as dual functional again $\calE$. Due to the lack of coerciveness and the absence of a bounded constraint set in $\calD$, it is not clear weather this functionals admits a minimizer. Nevertheless, it is possible to use $\calD$ to derive an error estimator directly, since the minimization of the dual functional is not necessary. \msa{Under these circumstances, }{} the solution of the primal problem would not enjoy the regularization approach and needs to be addressed via other techniques\msa{,}{} like \msa{e.g.}{} ADMM-method, see \cite{bib:BartelsMilicevicADMM, bib:BartelsMilicevicIterative}.\\
So far, we developed an abstract theoretical foundation for the a priori convergence as well as the a posteriori error estimation by the \msa{introduction}{} of primal-dual gap estimators. With these results and considerations at hand we close the theoretical discussion and draw our attention to the precise numerical utilization of \msa{our combined}{} findings to a selection of equilibrium problems.
\section{Numerical Application}
First, we introduce some common finite element spaces. Therefore, let a triangulation $\calT$ of an open, bounded domain $\Omega \subseteq \R^2$ be given. For given $k \in \N_0$ we define on a given \msa{triangle}{} $T \in \calT$ the space
\begin{eq*}
P_k(T) := \{u:T \to \R : u \text{ is a polynomial of degree} \leq k\}
\end{eq*}
and \msa{regarding the whole triangulation}{}
\begin{eq*}
P_k(\calT) := \left\{u \in L^{\msa{\infty}}(\Omega) : u|_T \in P_k(T) \text{ for all } T \in \calT \right\}.
\end{eq*}
The operator $\Pi_k :L^2(\Omega) \to P_k(\calT)$ denotes the metric projection (cf. \cite[Theorem 5.2]{bib:BrezisFunctionalAnalysis}) of the space $P_k(\calT) \subseteq L^2(\Omega)$ with respect to the $L^2$-norm. In the case $k = 0$ the projection reads as
\begin{eq*}
\Pi_0 u = \sum_{T \in \calT} \frac{1}{T} \int_T u \dx \cdot \mathbbm{1}_T.
\end{eq*}
For the discretization of objects in the Sobolev space $H^1(\Omega)$ and $H_0^1(\Omega)$ we define the spaces
\begin{eq*}
S_{(0)}^1(\calT) := \left\{ u \in H_{(0)}^1(\Omega) : u \in P_1(\calT) \right\}.
\end{eq*}
For the approximation of vector-valued objects\msa{, such}{} as \msa{gradients,}{} \msa{we}{} use the Raviart-Thomas finite element (see \cite{bib:RaviartThomas}) reading as
\begin{eq*}
\scalebox{0.975}{%
\ensuremath{%
RT_0(\calT) = \left\{ q \in H(\ddiv;\Omega) : \exists a_0 \in P_0(\calT;\R^d), a_1 \in P_0(\calT) :  q = a_0 + a_1 \cdot x \right\}.}}
\end{eq*}
This space can \msa{e.g. be used}{} for the mixed formulation of the Poisson problem
\begin{eq}[\label{eq:poisson}]
-\Delta y = f \text{ in } \Omega, \ y = 0 \text{ on } \partial \Omega.
\end{eq}
\msa{For its derivation}{}, one establishes a separate approximate for the gradient $p = \nabla y$ leading \msa{after partial integration}{} to the system:\\
Seek $p \in H(\ddiv,\Omega)$ and $y \in L^2(\Omega)$\msa{, such that}{}
\begin{eq}[\label{eq:poisoon:mixed}]
(p,q)_{L^2(\Omega;\R^d)} + (y,\ddiv q)_{L^2(\Omega)} &= 0 &&\text{for all } q \in H(\ddiv,\Omega)\msa{,}\\
(z,\ddiv p)_{L^2(\Omega)} &= -(f,z) &&\text{for all } z \in L^2(\Omega).
\end{eq}
One \msa{viable}{} discretization is the choice $P_0(\calT) \subseteq L^2(\Omega)$ and $RT_0(\calT) \subseteq H(\ddiv,\Omega)$ \msa{leading}{} to \msa{the following}{} stable conforming discretization \msa{(cf. \cite[Example 1.3.4]{bib:BoffiBrezziFortin})}:\\
Seek $p_\calT \in RT_0(\calT)$ and $y_\calT \in P_0(\calT)$ such that
\begin{eq}[\label{eq:poisson:mixed:discretized}]
(p_\calT,q_\calT)_{L^2(\Omega;\R^d)} + (y_\calT,\ddiv q_\calT)_{L^2(\Omega)} &= 0 &&\text{for all } q_\calT \in RT_0(\calT)\\
(z_\calT,\ddiv p_\calT)_{L^2(\Omega)} &= -(f,z_\calT) &&\text{for all } z_\calT \in P_0(\calT).
\end{eq}
By testing with $z_\calT = \mathbbm{1}_T$, denoting the characteristic function of a triangle $T \in \calT$, \msa{one deduces}{} $-\ddiv p_\calT = \Pi_0 f$. \msa{The mixed formulation}{} will be utilized in the upcoming subsection devoted to the detailed discussion of the error estimation explained by the example of the obstacle problem.
\subsection{Case Study: Obstacle Problem}
\msa{In the scope of this subsection we apply our combined findings to the obstacle problem, given in \refer{eq:obstacleproblem}.}{} At first, we notice, that in practice the imposition of a numerical discretization might come along with an inexact resolution of input data. \msa{In this situation one approximates}{} the primal \msa{and}{} dual functional, which leads to substitute functional $\calE_\gamma^\calT$ and $\calD_\gamma^\calT$\msa{, respectively}. \msa{This is investigated}{} for the instance of the obstacle problem. \msa{Here,}{} two input objects play a role, the force $f$ and the obstacle $\lowobs$, which might both \msa{need to be approximated}{} by $f_\calT \in \msa{P_0(\calT)}$ and $\lowobsh \in \msa{S^1(\calT)}$ with $\lowobsh|_{\partial \Omega} < 0$. Then, we formulate the functionals
\begin{eq*}
\calE_\gamma^\calT(y_\calT) := \frac{1}{2}\|\nabla y_\calT\|_{L^2(\Omega;\R^d)}^2 - (f_\calT, i_{L^2(\Omega)}y_{\msa{\calT}})_{L^2(\Omega)} + \frac{\gamma}{2}\int_\Omega (\lowobsh - y_\calT)^{2+} \dx
\end{eq*}
for $y_\calT \in S_0^1(\calT)$ and
\begin{eq*}
\calD_\gamma^\calT(p_\calT,z_\calT) &:= \frac{1}{2}\|p_\calT\|_{L^2(\Omega;\R^d)}^2 + (z_\calT,\lowobsh)_{L^2(\Omega)} + \frac{1}{2\gamma}\|z_\calT\|_{L^2(\Omega)}^2\\
&+ I_{L^2_-(\Omega)}(z_\calT) + I_{\{0\}}(-\ddiv p_\calT + z_\calT - f_\calT)
\end{eq*}
for $p_\calT \in RT_0(\calT)$ and $z_\calT \in P_0(\calT)$. For the calculation of the estimator we define \msa{$\eta_{\gamma,\calT}^2(y_\calT,p_\calT,z_\calT) = \calE_\gamma^\calT(y_\calT) + \calD_\gamma^\calT(p_\calT,z_\calT)$}{} and decompose
\msa{\begin{eq}[\label{eq:error:decomp}]
\eta_\gamma^2(y_\calT,p,z) &= \calE_\gamma(y_\calT) + \calD_\gamma(p,z)\\
&= \eta_{\gamma,\calT}^2(y_\calT,p_\calT,z_\calT) + \left(\calE_\gamma(u_\calT) - \calE_\gamma^\calT(y_\calT)\right)\\
&+ \left(\calD_\gamma(p,z) - \calD_\gamma^\calT(p_\calT,z_\calT)\right).
\end{eq}}{}%
\msa{Here, $p \in H(\ddiv,\Omega)$ and $z \in L^2(\Omega)$ can be chosen arbitrary in the above calculation.}{} In practice\msa{,}{} we \msa{only choose $p_\calT \in RT_0(\calT)$ and $z_\calT \in P_0(\calT)$}{}, and calculate $\eta_{\gamma,\calT}^2$. \msa{It is important to note here, that for the primal functional we just plugged in the same argument. In contrast, the dual functionals include an indicator function containing an operator equation, which needs to be fulfilled \emph{exactly}. Thus, having discretized data and being limited to a discretized subspace necessitates the presence of the functions $p$ and $z$.}{}\\
\msa{In \refer{eq:error:decomp} t}{}he second and third term contributes as well to the error and hence needs to be estimated\msa{. This}{} is investigated in the following Theorem.
\begin{thm}\label{thm:oscillation}
Let arbitrary $y_\calT \in S_0^1(\calT)$, $p_\calT \in RT_0(\calT)$ and $z_\calT \in P_0(\calT)$ be given. Then, there exist $p \in H(\ddiv,\Omega)$ and $z \in L^2(\Omega)$\msa{,}{} such that
\begin{eq*}
\calE_\gamma(\msa{y}_\calT) - \calE_\gamma^\calT(\msa{y}_\calT) &\lesssim \|f - f_\calT\|_{L^2(\Omega)}\|y_\calT\|_{L^2(\Omega)}\\
&+ \gamma \|(\lowobs - y_\calT)^+\|_{L^2(\Omega)} \|\lowobs - \lowobsh\|_{L^2(\Omega)} \text{ and}\\
\calD_\gamma(p,z) - \calD_\gamma^\calT(p_\calT,z_\calT) &\lesssim \left( \|f- f_\calT\|_{L^2(\Omega)} + \|p_\calT\|_{L^2(\Omega;\R^d)} \right)\|f - f_\calT\|_{L^2(\Omega)}\\
&+ \|z_\calT\|_{L^2(\Omega)}\|\lowobs - \lowobsh\|_{L^2(\Omega)}
\end{eq*}
\msa{hold true.}
\end{thm}
\begin{proof}
First, we obtain
\begin{eq*}
(f - f_\calT,y_\calT)_{L^2(\Omega)} \lesssim \|f - f_\calT\|_{L^2(\Omega)} \|y_\calT\|_{L^2(\Omega)}.
\end{eq*}
Using the inequality $a^2 - b^2 \leq 2a(a-b)$ for $a,b \in \R$ we obtain
\begin{eq*}
\frac{\gamma}{2}\int_\Omega (\lowobs &- y_\calT)^{2+} \dx - \frac{\gamma}{2}\int_\Omega (\lowobsh - y_\calT)^{2+} \dx\\
&\leq \gamma \|(\lowobs - y_\calT)^+\|_{L^2(\Omega)} \|(\lowobs - y_\calT)^+ - (\lowobsh - y_\calT)^+\|_{L^2(\Omega)}\\
&\leq \gamma \|(\lowobs - y_\calT)^+\|_{L^2(\Omega)} \|\lowobs - \lowobsh\|_{L^2(\Omega)}.
\end{eq*}
The combination of these estimates leads to the assertion for the primal functionals $\calE_\gamma, \calE_\gamma^\calT$.
For the second inequality we take $z = z_\calT$ and solve the equation
\begin{eq}[\tag{$*$}\label{eq:osc:tmp}]
- \Delta \xi = f - f_\calT \text{ in } \Omega,\quad \xi = 0 \text{ on } \partial \Omega
\end{eq}
and define $p = p_\calT + \nabla \xi$. Using $a^2 - b^2 \leq 2 a (a-b)$ and $z = z_\calT$  we obtain
\begin{eq*}
\calD_\gamma(p,z) - \calD_\gamma^\calT(p_\calT,z_\calT) &= \frac{1}{2}\|p\|_{L^2(\Omega;\R^d)}^2 - \frac{1}{2}\|p_\calT\|_{L^2(\Omega;\R^d)}^2\\
&+ (z_\calT, \lowobs - \lowobsh)_{L^2(\Omega)}\\
&\leq \|p\|_{L^2(\Omega;\R^d)} \|p - p_\calT\|_{L^2(\Omega;\R^d)}\\
&+ \|z_\calT\|_{L^2(\Omega)} \|\lowobs - \lowobsh\|_{L^2(\Omega)}.
\end{eq*}
Testing \refer{eq:osc:tmp} with $\xi$ itself leads to
\begin{eq*}
\|\nabla \xi\|_{L^2(\Omega;\R^d)}^2 &= (f - f_\calT, \xi)_{L^2(\Omega)} \leq \|f - f_\calT\|_{L^2(\Omega)} \|\xi\|_{L^2(\Omega)}\\
&\lesssim \|f- f_\calT\|_{L^2(\Omega)}\|\nabla \xi\|_{L^2(\Omega;\R^d)}.
\end{eq*}
Subsequently, we get
\begin{eq*}
\calD_\gamma(p,z) - \calD_\gamma^\calT(p_\calT,z_\calT) &\lesssim \left(\|p_\calT\|_{L^2(\Omega;\R^d)} +  \|\nabla \xi\|_{L^2(\Omega;\R^d)}\right)\|\nabla \xi\|_{L^2(\Omega;\R^d)}\\
&+ \|z_\calT\|_{L^2(\Omega)}\|\lowobs - \lowobsh\|_{L^2(\Omega)},
\end{eq*}
which leads to the remaining \msa{assertion}.
\end{proof}
The proposed estimate suggests, that in the case of an approximated obstacle, the boundedness of the Moreau--Yosida regularization $\frac{\gamma}{2}\int_\Omega (\lowobs - y_\calT)^{2+} \dx$ leads to the \msa{requirement}{} $\sqrt{\gamma}\|\lowobs - \lowobsh\|_{L^2(\Omega)} \to 0$ \msa{to guarantee overall convergence to zero}. The choice of $y_\calT$ is of course based on the solution of the equilibrium problem associated to $\calE_\gamma^\calT$. This can be achieved via a semi-smooth Newton method \cite{bib:ChenNashedZuhair, bib:HintItoKunisch} applied to the first order system of the discretized problem. To obtain a sharp estimate, one might be tempted to minimize the dual functional $\calD_\gamma^\calT$ as well, as it has been done for related problems (see e.g. \cite{bib:BartelsMilicevicADMM,bib:BartelsMilicevicPrimalDualGap}).
However, we would like to avoid this situation for several reasons: In contrast to optimization problems, the dual problem incorporates the equilibrium itself. This imposes a hierarchy between these problems \msa{needing the primal problem to be solved first and the dual problem afterwards}{}. Thus, \msa{both problems}{} can in general not be solved in parallel. Moreover, looking at the shape of the dual functional in the instance of the obstacle problem, it is difficult to apply the same strategy based on a semi-smooth Newton method as for the primal problem, since an obstacle-type constraint in $L^2(\Omega)$ is imposed, but no additional regularity can be exploited.
\msa{Thus, other numerical methods}{}, e.g. ADMM \msa{(}as successfully applied in \cite{bib:BartelsMilicevicADMM}\msa{)}{} \msa{are of interest}. As observed therein, the solution of the dual problem took more iterations in comparison to the primal one\msa{,}{} and \msa{the solution of}{} the adaptive case took more iterations compared to uniform refinement. Hence, the minimization of the dual problem for the primal-dual gap estimate might be computationally expensive.\\
For the sake of completeness we remark here, that the problem of constructing a suitable point to be plugged in the dual formulation has been addressed in \cite{bib:BartelsNonconformingPrimalDualGap}. Therein, the relationship between the Crouzeix--Raviart and Raviart--Thomas finite element solution (as in \cite{bib:Marini}) are generalized and used to derive estimators via the first order system.\\ 
However, we go here for a simpler approach instead\msa{:}{} As we are allowed to plug in any choice of feasible point to obtain a reliable error estimator it might already be sufficient to plug in a point, that is based on the calculation of the approximated equilibrium $y_\calT$. To do so, we recall the first order system for a minimizer of the Moreau--Yosida regularized obstacle problem
\begin{eq*}
-\Delta y - \gamma (\lowobs-y)^+ = f \text{ in }\Omega,\quad y = 0 \text{ on } \partial \Omega.
\end{eq*}
For the dual problem the relations $p = \nabla y$ and $z = -\gamma(\lowobs-y)^+$ hold. As discretized analogue we choose $z_\calT = -\gamma \Pi_0(\lowobsh - y_\calT)^+ \in P_0(\calT)$ and solve the mixed formulation of the Poisson equation \msa{\refer{eq:poisson:mixed:discretized}}{} with right hand side $f_\calT - \gamma (\lowobsh - y_\calT)^+$ based on the Raviart-Thomas finite element. This leads to the following linear system:\\
Seek $p_\calT \in RT_0(\calT)$ and $v_\calT \in P_0(\calT)$ such that
\begin{eq*}
(p_\calT,q_\calT)_{L^2(\Omega,\R^d)} + (v_\calT,\ddiv q_\calT)_{L^2(\Omega)} &= 0 \text{ for all } q_\calT \in RT_{0}(\calT),\\
(\ddiv p_\calT, w_\calT)_{L^2(\Omega)} + (f_\calT + \gamma(\lowobsh - y_\calT)^+, w_\calT)_{L^2(\Omega)} &= 0  \text{ for all } w_\calT \in P_0(\calT).
\end{eq*}
Then\msa{,}{} $-\ddiv p_\calT = f_\calT + \gamma \Pi_0 (\lowobsh - y_\calT)^+$ holds \msa{as $f_\calT \in P_0(\calT)$ is assumed}. Together with partial integration this observation leads to
\begin{eq*}
(f_\calT,y_\calT)_{L^2(\Omega)} &= (-\ddiv p_\calT \msa{-}{} \gamma \Pi_0 (\lowobsh - y_\calT)^+,y_\calT)_{L^2(\Omega)}\\
&= (p_\calT, \nabla y_\calT)_{L^2(\Omega;\R^d)} + (z_\calT,y_\calT)_{L^2(\Omega)},
\end{eq*}
and thus to the estimator 
\begin{eq*}
\eta_{\gamma,\calT}^2&(y_\calT,p_\calT,z_\calT) = \calE_\gamma^\calT(y_\calT) + \calD_\gamma^\calT(p_\calT,z_\calT)\\
&= \frac{1}{2}\|\nabla y_\calT\|_{L^2(\Omega;\R^d)}^2 - (f_\calT,y_\calT)_{L^2(\Omega)} + \frac{\gamma}{2}\int_\Omega (\lowobsh - y_\calT)^{2+} \dx\\
&+ \frac{1}{2}\|p_\calT\|_{L^2(\Omega;\R^d)}^2 + \frac{\gamma}{2}\|\Pi_0(\lowobsh - y_\calT)^+\|_{L^2(\Omega)}^2 + (z_\calT,\lowobsh)_{L^2(\Omega)} \\
&= \frac{1}{2}\|\nabla y_\calT\|_{L^2(\Omega;\R^d)}^2 + \frac{1}{2}\|p_\calT\|_{L^2(\Omega;\R^d)}^2 - (\nabla y_\calT, p_\calT)_{L^2(\Omega;\R^d)}\\
&+ \frac{\gamma}{2}\|(\lowobsh - y_\calT)^+\|_{L^2(\Omega)}^2 + \frac{\gamma}{2}\|\Pi_0 (\lowobsh - y_\calT)^+\|_{L^2(\Omega)}^2\\
&- \gamma(\Pi_0(\lowobsh - y_\calT)^+,\lowobsh - y_\calT)_{L^2(\Omega)}\\
&= \frac{1}{2}\|p_\calT - \nabla y_\calT\|_{L^2(\Omega;\R^d)}^2 + \frac{\gamma}{2}\|( \lowobsh - y_\calT)^+ - \Pi_0(\lowobsh - y_\calT)^+\|_{L^2(\Omega)}^2\\
&+ \gamma(\Pi_0(\lowobsh - y_\calT)^+,(y_\calT - \lowobsh)^+)_{L^2(\Omega)}.
\end{eq*}
This decomposition can be interpreted as follows: The first term compares the computed approximation with another one with respect to a higher order \msa{element}, as the variable $p_\calT$ is piecewise linear whereas $\nabla y_\calT$ is piecewise constant.
\msa{Implicitly, this}{} term encodes a residual quantity as well as an edge part. The second term compares the surplus of the state with respect to the lower obstacle with its piecewise constant approximation and hence builds a part solely \msa{focussing}{} the inactive set of the computed state. The third term is structurally similar, but considers those triangles containing points of the inactive as well as active sets and is hence devoted to their interface region. 
For the last term we obtain
\begin{eq*}
\gamma (\Pi_0&(\lowobsh - y_\calT)^+,(y_\calT - \lowobsh)^+)_{L^2(\Omega)}\\
&= \gamma (\Pi_0(\lowobsh - y_\calT)^+ - (\lowobsh - y_\calT)^+, (y_\calT - \lowobsh)^+ )_{L^2(\Omega)}\\
&= \gamma (\Pi_0(\lowobsh - y_\calT)^+ - (\lowobsh - y_\calT)^+, (y_\calT - \lowobsh)^+ - \Pi_0(y_\calT - \lowobsh)^+)\msa{.}{}
\end{eq*}
The choice of $p_\calT,z_\calT$ is in principle heuristical and \msa{likely}{} suboptimal\msa{, but still}{} based on the calculated approximate solution and \msa{using the first order system and is hence expected to improve}{} for fine meshes.
\subsection{Adaptive Refinement}
In our approach we \msa{use}{} regularization and discretization schemes simultaneously. In the light of the Diagonal theorem in \refer{thm:diag} the selection of a subsequence is required. This \msa{was in its proof}{} based on balancing of the approximation error of the equilibrium, which is encoded in the error estimators as well as the approximation errors discussed in \refer{thm:oscillation}.

From a practical viewpoint we \msa{require}{} a strategy to decide when to perform another mesh refinement and when to update the regularization parameter instead.
\msa{Besides the above mentioned Diagonal theorem, other performance reasons play a role.}{}
Relating to the obstacle problem the solution \msa{might}{} have areas, where the values of the solution and the obstacle coincide. For \msa{them}, early updates of the regularization parameters are relevant, because otherwise many degrees of freedoms might be invested by the algorithm to address areas relevant for small $\gamma$, that become irrelevant for the later stages of the method. An aggressive refinement strategy \msa{might lead}{} for late $\gamma$-updates to \msa{convergence}{} issues with the applied solution method, especially with the basin of convergence for the semi-smooth Newton method.\\
A successfully applied strategy for the update of $\gamma$ is the formulation of a value functional (see \cite{bib:HintKunischPath}). \msa{This has been successfully applied to}{} optimization problems\msa{, where the value function is}{} defined as the minimum of the objective depending on $\gamma$.
\msa{In the scope of this work however,}{} we exploit \refer{eq:aposteriori} and use the error estimator itself as value functional. The idea is to control the increase \msa{of}{} the error estimator induced by the update of the regularization parameter. Therefore, we \msa{demand}{} the estimator for the new parameter $\gamma_{n+1}$ to differ from the one on the old level $\gamma_n$ (on the same mesh) \msa{at most}{} by a factor of $(1+c_\gamma)$ with a freely chosen $c_\gamma > 0$. To achieve this, one performs a \msa{(formal)}{} linearization via
\begin{eq*}
\eta_{\gamma_{n+1},\ell}^2 \approx \eta_{\gamma_n,\ell}^2 + (\gamma_{n + 1} - \gamma_n)\frac{d}{d\gamma}\eta_{\gamma_n,\ell}^2.
\end{eq*}
Here \msa{and in the following}{}, we write for a sequence of triangulations $(\calT_\ell)_{\ell \in \N}$ rather the index $\ell$ instead of $\calT_\ell$. For the calculation of the derivative we use \msa{\refer{lem:duality:my}}{} and exploit the structure of the theoretically optimal error estimator reading as
\begin{eq*}
\eta_{\gamma,\ell,\mathrm{opt}}^2 := \inf_{p_\ell \in \calP_\ell}\left( \calE_\gamma^\ell(u_\ell,u_\ell) + \calD_\gamma^\ell(p_\ell,u_\ell)\right)
\end{eq*}
with the functional being of the type $f(x) + \gamma \pi(x) + \frac{1}{\gamma}\rho(x)$ \msa{for respective functionals $f,\rho,\pi$}. Under some mild assumptions (see \refer{lem:apx:derivative}) \msa{imposed on the latter,}{} one can then derive an expression for the derivative with respect to the regularization parameter of the optimal choice. Together with the \msa{condition}{} $\eta_{\gamma_{n+1},\ell}^2 \leq (1+c_\gamma)\eta_{\gamma_n,\ell}^2$ one deduces
\begin{eq*}
\gamma_{n+1} \lessapprox \gamma_n + c_\gamma\frac{\eta_{\gamma_n,\ell}^2}{\frac{d}{d\gamma}\eta_{\gamma_n,\ell}^2},
\end{eq*}
provided the derivative is positive. Additionally, to this approximation we safeguard the update from below with a parameter $\gamma_{\mathrm{min}\,\mathrm{update}}$ \msa{to avoid a rapid sequence of small updates at the beginning of the algorithm}{}. In the case of a non-positive derivative, we just multiply the regularization parameter with the constant factor $(1+c_\gamma)$ instead. To address these considerations we formulate the following Algorithm.\\
\begin{algorithm}[H]
    \SetAlgoLined
    \KwData{Penalty parameter $0 < \gamma_0$, $0 < c_\gamma, 0 < c_\eta < 1$, start triangulation $\calT_0$, Doerfler parameter $\theta \in (0,1)$, $\gamma_{\msa{\mathrm{min}\,\mathrm{update}}}, \gamma_{\mathrm{max}}, \text{\tt nrdof}_{\mathrm{max}}> 0$}
    \For{$n = 0, 1, \dots$}{
        \For{$\ell = 0, 1, \dots$}{
        Solve the equilibrium problem for $\calE^\ell_{\gamma_n}$\;     
        \If{$\gamma_n \geq \gamma_{\mathrm{max}}$ or $\texttt{\emph{nrdof}} \geq \texttt{\emph{nrdof}}_{\mathrm{max}}$}{\textbf{return}\;}
        Calculate the error estimator $\eta_{\gamma_n,\ell}$\;
        	\If{$\ell = 1$ and $n = 1$}{
        		Set $\eta_{\mathrm{ref}} = \eta_{\gamma_n,\ell}$
        	}
        	\eIf{$\eta_{\gamma_n,\ell} \leq c_\eta \eta_{\mathrm{ref}}$}{
        		$\eta_{\mathrm{ref}} = \eta_{\gamma_n,\ell}$\;
        		\eIf{$\frac{d}{d\gamma}\eta_{\gamma_n, \ell}^2 > 0$}{$\gamma_{n+1} = \gamma_n + \max\left(\gamma_{\mathrm{min}\, \mathrm{update}}, c_{\gamma}\frac{\eta_{\gamma_n,\ell}^2}{\frac{d}{d\gamma} \eta^2_{\gamma_n,\ell}}\right)$\;}{$\gamma_{n+1} = (1 + c_\gamma) \gamma_n$}
        		Set $\calT_0 = \calT_\ell$\;
        		\textbf{break}
        	}{$\calT_{\ell + 1} = \operatorname{\tt doerfler\_refine}(\theta,\eta_{\gamma_n,\ell}^2,\calT_{\ell})$}
        }
    }
\caption{Joint mesh and $\gamma$-update strategy}\label{alg:afem}
\end{algorithm}
Here, we use Doerfler marking (cf. \cite{bib:Doerfler}) as refinement strategy. The \msa{parameter}{} $c_\eta$ controls the reduction of the estimator to be achieved by the mesh refinement and is supposed to be a number between zero and one (in our experiments we used $\frac{1}{2}$). The comparison relates to the lowest estimator value of the old mesh. This has been done to guarantee the values of the estimators \msa{to}{} go to zero in total.\\
In the following, we \msa{apply}{} our combined findings to a selction of quasi-variational inequalities (abbr.: QVIs) in Subsections \ref{ssec:qvi_thermo} and \ref{ssec:qvi_membrane}. Those \msa{applications}{} cover obstacle-type constraints and will utilize as well our findings for the obstacle problem.\\
Our notation might differ in the following from the one used so far in this article, since we want to follow the notation given in the original works containing these applications as closely as possible.
\subsection{A QVI with Application to Thermoforming}\label{ssec:qvi_thermo}
\begin{figure}[t]
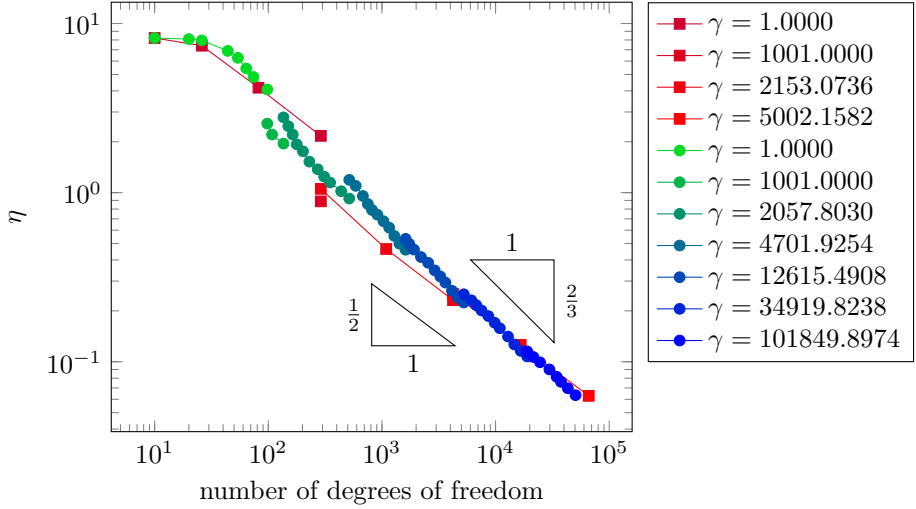

\newcommand\datapath{plots/qvi_thermo/data/square}%
\includestandalone{plots/qvi_thermo/data/square/qvi_thermo_est_final}%

\caption{Values of the error estimators for different values of $\gamma$ (see color code on the right). The adaptive refinement and update strategy (circles) is compared with the uniform refinement (squares).}\label{fig:qvi_thermo:afem}
\end{figure}

\begin{figure}
\includegraphics[trim = {890, 250, 890, 250}, clip, height = 0.217\textheight]{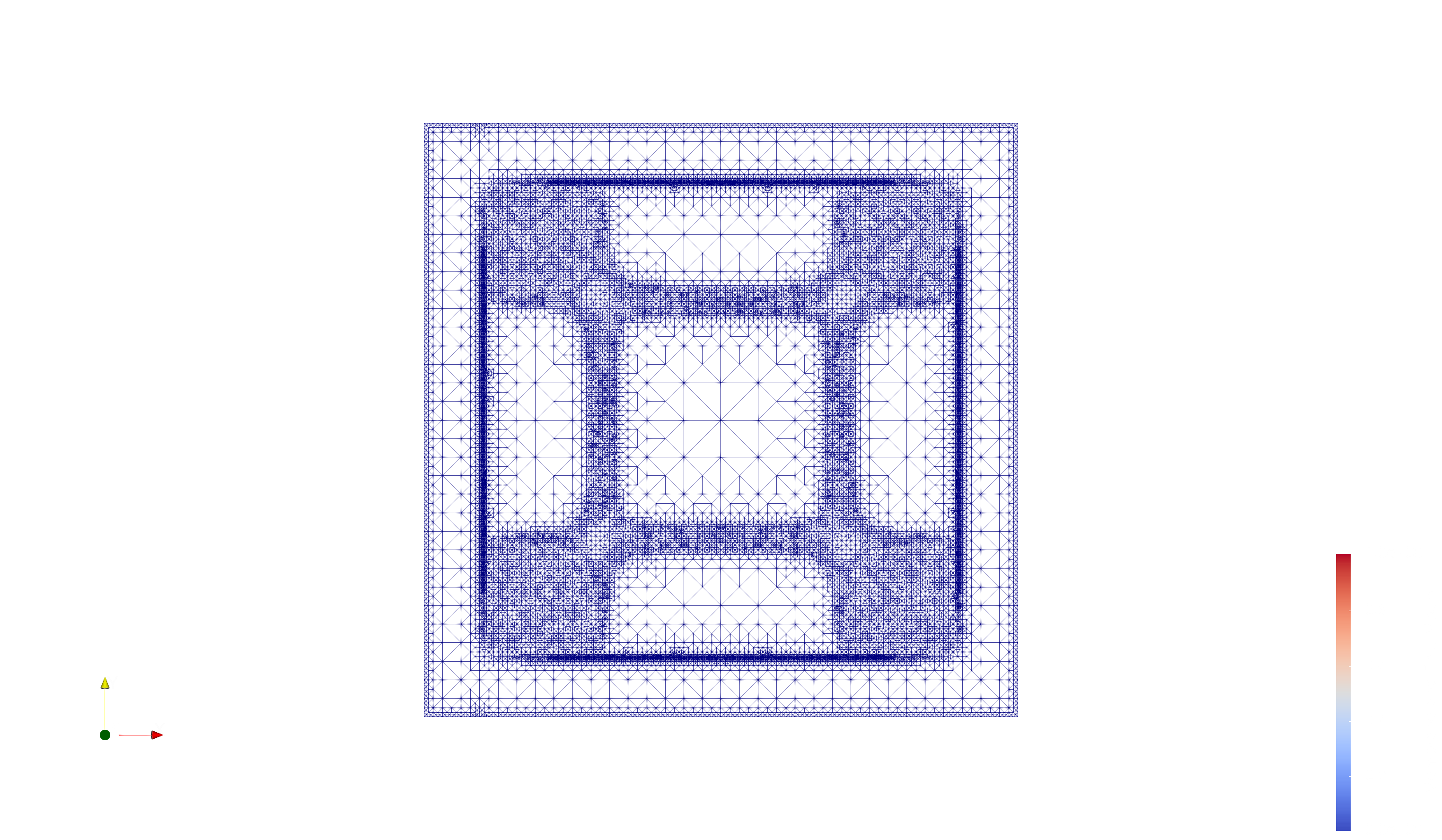}
\includegraphics[trim = {650, 200, 300, 450}, clip, height = 0.217\textheight]{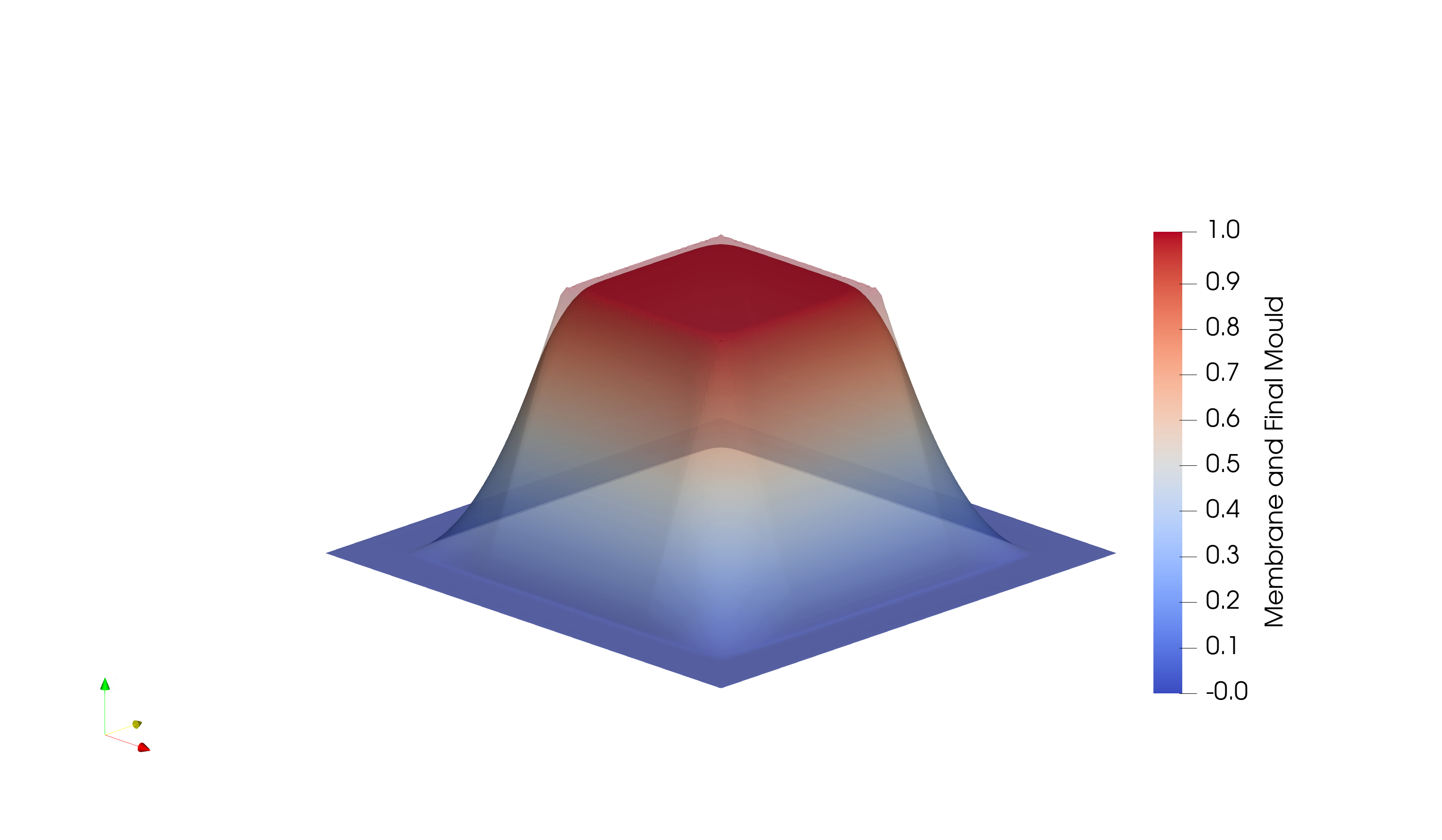}
\caption{Left: Adaptively refined mesh. Right: Plot of the final mould $\Phi = \Phi_0 + LT$ and the membrane.}\label{fig:qvi_thermo:mesh_plot}
\end{figure}

As first example, we \msa{apply}{} our results to the numerical treatment of the following quasi-variational inequality, discussed in \cite[Section 6]{bib:AlphonseHintRautenbergDirectDiff}\msa{:}{}\\
Seek $y \in H_0^1(\Omega)$ and $T \in H^1(\Omega)$ such that
\begin{sys}[\label{sys:qvi_thermo}]
f &\in -\Delta u  + N_{C(u)}(u),\label{eq:sys:qvi_thermo:vi}\\
k T - \Delta T &= g(\Phi(u) - u) \text{ in } \Omega,\label{eq:sys:qvi_thermo:pde}\\
\frac{\partial T}{\partial \nu} &= 0 \text{ on } \partial \Omega\label{eq:sys:qvi_thermo:nbc}
\end{sys}
holds \msa{where}{} $C(u) = \{ v \in H_0^1(\Omega) : v \leq \Phi(u) \}$ \msa{for an operator $\Phi:U \to U$ with $\Phi(u) = \Phi_0 + LT$}{} and \msa{$g:\R\to \R$ being}{} a globally Lipschitz continuous\msa{,}{} twice continuously differentiable \msa{and}{} decreasing function. Using \msa{the definition of $\Phi$}{} we can write $C_u(T):= \{v \in H_0^1(\Omega) : v \leq \Phi_0 + LT\}$ with $L \in \calL(H^1(\Omega),H_0^1(\Omega))$ \msa{being the multiplication with an $L^\infty$-function. For the details on the parameter choices we refer again to \cite{bib:AlphonseHintRautenbergDirectDiff}. Therein,}{} the existence has been derived using the Birkhoff--Tartar theorem (cf. \cite[Chapter 2.5]{bib:MoscoQVI}).\\
\msa{For the embedding into our framework, we interpret \refer{sys:qvi_thermo}}{} as the first order system of the equilibrium problem in the sense of \refer{eq:equi} using the functional
\begin{eq*}
\calE(v,S; u,T) :&= \frac{1}{2}\|\nabla v\|_{L^2(\Omega;\R^d)}^2 - (f,v)_{\msa{L^2(\Omega)}} + I_{C_u(T)}(v)\\
&+ \frac{k}{2}\|S\|_{L^2(\Omega)}^2 + \frac{1}{2}\|\nabla S\|_{L^2(\Omega;\R^d)}^2 - (g(\Phi_0 + LT-u),S)_{L^2(\Omega)}.
\end{eq*}
Just as for the obstacle problem, we reinterpret the obstacle condition in $L^2(\Omega)$. Returning to the discussion after \refer{lem:my:gamma} one can clearly see, that the conditions therein are fulfilled as $L$ is a bounded, linear operator. To use the developed framework, we \msa{set}{}
\begin{eq}[\label{eq:qvi_thermo:hfunctional}]
&H: L^2(\Omega) \times H_0^1(\Omega) \times H^1(\Omega) \to [0,+\infty] \text{ defined by}\\
&H(z;u,T):= I_{L_+^2(\Omega)}(\Phi_0 + LT - z).
\end{eq}
Then, we obtain as Moreau--Yosida regularized functional
\begin{eq*}
\calE_\gamma(v,S;u,T) &:= \frac{1}{2}\|\nabla v\|_{L^2(\Omega;\R^d)}^2 - (f,v)_{\msa{L^2(\Omega)}} + \frac{\gamma}{2}\int_\Omega (v - \Phi_0 - LT)^{2+} \dx\\
&+ \frac{k}{2}\|S\|_{L^2(\Omega)}^2 + \frac{1}{2}\|\nabla S\|_{L^2(\Omega;\R^d)}^2 - (g(\Phi_0 + LT-u),S)_{L^2(\Omega)}.
\end{eq*}
The latter leads to the first order condition
\begin{sys}[\label{sys:qvi_thermo:pen}]
-\Delta u + \gamma (u - \Phi_0 - LT)^+ &= f &&\text{in } \Omega,\label{eq:sys:qvi_thermo:pen:vi}\\
u &= 0 &&\text{on } \partial \Omega,\label{eq:sys:qvi_thermo:pen:dbc}\\
kT - \Delta T &= g(\Phi_0 + LT - u) &&\text{in } \Omega,\label{eq:sys:qvi_thermo:pen:pde}\\
\frac{\partial T}{\partial \nu} &= 0 &&\text{on } \partial \Omega.\label{eq:sys:qvi_thermo:pen:nbc}
\end{sys}
With the reasoning in \cite[Section 3.1 and Section 6]{bib:AlphonseHintRautenbergDirectDiff} it is straightforward to show the existence of a bounded sequence of solutions of the regularized system \msa{as $\gamma \to \infty$}{}.\\
\msa{To}{} apply \refer{lem:duality:my}\msa{, take}{} the spaces $U := H_0^1(\Omega) \times H^1(\Omega)$ and $Y:= L^2(\Omega;\R^d) \times L^2(\Omega) \times L^2(\Omega;\R^d) \times L^2(\Omega) \times L^2(\Omega)$ as well as $Z = L^2(\Omega)$\msa{. For the functionals we set}{} $F(v,S;u,T) := 0$ and
\begin{eq*}
G(p,z_v,q,z_{S,1},z_{S,2}) :&= \frac{1}{2}\|p\|_{L^2(\Omega;\R^d)}^2 -(f,z_v)_{L^2(\Omega)} + \frac{1}{2}\|q\|_{L^2(\Omega;\R^d)}^2\\
&+ \frac{k}{2}\|z_{S,1}\|_{L^2(\Omega)}^2 - (g(\Phi_0 + LT - u),z_{S,2})_{L^2(\Omega)}
\end{eq*}
\msa{with $H$ as in \refer{eq:qvi_thermo:hfunctional}.}{} \msa{As operators we choose}{} $A(v,S):= (\nabla v, v, \nabla S, S, S)$ \msa{and}{} $B(v,S) := i_{L^2(\Omega)}v$. Using Calculation \ref{calc:conjugate} the dual functional reads after slight simplifications as
\begin{eq}
\calD_\gamma&(p,s,q,z,;u,T) = \frac{1}{2}\|p\|_{L^2(\Omega;\R^d)}^2 + \frac{1}{2}\|q\|_{L^2(\Omega;\R^d)}^2 + \frac{1}{2k}\|s\|_{L^2(\Omega)}^2\\
&+ (\Phi_0 + LT,z)_{L^2(\Omega)} + \frac{1}{2\gamma}\|z\|_{L^2(\Omega)}^2 + I_{\{0\}}(-\ddiv p + z^* - f)\\
&+ I_{\{0\}}(-\ddiv q + s^* - g(\Phi_0 + LT - u)).
\end{eq}
Just as for the obstacle problem, we need a discretization of the input data. Since in \cite[Section 6]{bib:AlphonseHintRautenbergDirectDiff}, $f$ is assumed to be constant, we just need to employ a discretization on the initial mould $\Phi_0$ via a function $\Phi_{0,\calT}$. Using the reasoning in \refer{thm:oscillation}, we obtain with
\begin{eq*}
\calE^\calT_\gamma(v_\calT,S_\calT;u_\calT,T_\calT) :&= \frac{1}{2}\|\nabla v_\calT\|_{L^2(\Omega;\R^d)}^2 - (f,v_\calT)_{L^2(\Omega)} + \frac{k}{2}\|S_\calT\|_{L^2(\Omega)}^2\\
&+ \frac{1}{2}\|\nabla S_\calT\|_{L^2(\Omega;\R^d)}^2 + \frac{\gamma}{2}\int_\Omega (v_\calT - \Phi_{0,\calT} - LT_\calT)^{2+} \dx\\
& - (g(\Phi_{0,\calT} + LT_\calT-u_\calT),S_\calT)_{L^2(\Omega)}
\end{eq*}
that
\begin{eq*}
\calE_\gamma(v_\calT,S_\calT;u_\calT,T_\calT) &- \calE_\gamma^\calT(v_\calT,S_\calT;u_\calT,T_\calT)\\
&\lesssim \gamma \|(v_\calT - \Phi_0 - LT_\calT)^+\|_{L^2(\Omega)} \|\Phi_0 - \Phi_{0,\calT}\|_{L^2(\Omega)}\\
&+ \Lip{g}\|\Phi_0 - \Phi_{0,\calT}\|_{L^2(\Omega)} \|S_\calT\|_{L^2(\Omega)}
\end{eq*}
and with
\begin{eq*}
&\calD_\gamma^\calT(p_\calT,s_\calT,q_\calT,z_\calT;u_\calT,T_\calT) = \frac{1}{2}\|p_\calT\|_{L^2(\Omega;\R^d)}^2 + \frac{1}{2}\|q_\calT\|_{L^2(\Omega;\R^d)}^2\\
&+ \frac{1}{2k}\|s_\calT\|_{L^2(\Omega)}^2 + (\Phi_{0,\calT} + LT_\calT,z_\calT)_{L^2(\Omega)} + \frac{1}{2\gamma}\|z_\calT\|_{L^2(\Omega)}^2\\
&+ I_{\{0\}}(-\ddiv p_\calT + z_\calT - f) + I_{\{0\}}(-\ddiv q_\calT + s_\calT - g(\Phi_{0,\calT} + LT_\calT - u_\calT))
\end{eq*}
that for $p_\calT,q_\calT \in RT_0(\calT)$ and $z_\calT,s_\calT \in P_0(\msa{\calT})$ there exist $p,q \in L^2(\Omega;\R^d)$ and $z,s \in L^2(\Omega)$ such that
\begin{eq*}
&\calD_\gamma(p,s,q,z;u_\calT,T_\calT) - \calD_\gamma^\calT(p_\calT,s_\calT,q_\calT,z_\calT;u_\calT,T_\calT)\\
&\lesssim \Lip{g}\left(\Lip{g}\|\Phi_{0,\calT} - \Phi_0\|_{L^2(\Omega)} + \|q_\calT\|_{L^2(\Omega;\R^d)}\right)\|\Phi_{0,\calT} - \Phi_0\|_{L^2(\Omega)}\\
&+ \|z_\calT\|_{L^2(\Omega)}\|\Phi_{0,\calT} - \Phi_0\|_{L^2(\Omega)},
\end{eq*}
where $z = z_\calT, s = s_\calT$ as well as $p = p_\calT$ and $q$ as in the proof of \refer{thm:oscillation} have been chosen.\\
\msa{Next, we}{} derive the error estimator.
\msa{For that, we set $z_\calT = \gamma \Pi_0(u_\calT - L T_\calT - \Phi_{0,\calT})^+$ and $s_\calT = k \Pi_0 T_\calT$ and take $p_\calT,q_\calT \in RT_0(\calT)$ to be the solutions of the mixed formulation of the Poisson equation with right hand sides $f - z_\calT$ and $g(\Phi_{0,\calT} + LT_\calT - u_\calT) - s_\calT$  respectively. Using partial integration we obtain with $-\ddiv p_\calT + z_\calT = f$ the identity}{}
\begin{eq*}
(f,u_\calT)_{L^2(\Omega)} = (z_\calT,u_\calT)_{L^2(\Omega)} + (p_\calT, \nabla u_\calT)_{L^2(\Omega;\R^d)}
\end{eq*}
and with $-\ddiv q_\calT + s_\calT = g(\Phi_{0,\calT} + L T_\calT - u_\calT)$ \msa{the relation}
\begin{eq*}
(g(u_\calT - \Phi_{0,\calT} - L T_\calT),T_\calT)_{L^2(\Omega)} = (s_\calT,u_\calT)_{L^2(\Omega)} + (q_\calT, \nabla T_\calT)_{L^2(\Omega;\R^d)}\msa{.}
\end{eq*}
\msa{Thus, the estimator reads as}{}
\begin{eq*}
\eta_{\gamma,\calT}^2 :&= \calE_\gamma^\calT(u_\calT,T_\calT;u_\calT,T_\calT) + \calD_\gamma^\calT(p_\calT,s_\calT,q_\calT,z_\calT;u_\calT,T_\calT)\\
&= \frac{1}{2}\|\nabla u_\calT - p_\calT\|_{L^2(\Omega;\R^d)}^2 + \frac{1}{2}\|\nabla T_\calT - q_\calT\|_{L^2(\Omega;\R^d)}^2\\
&+ \frac{1}{2k}\|k T_\calT - s_\calT\|_{L^2(\Omega)}^2 + \frac{\gamma}{2}\int_\Omega (u_\calT - \Phi_{0,\calT} - LT_\calT)^{2+} \dx\\
&+ \frac{1}{2\gamma}\|z_\calT\|_{L^2(\Omega)}^2 - (z_\calT, u_\calT - \Phi_{0,\calT} - L T_\calT)_{L^2(\Omega)},
\end{eq*}
where we suppressed the dependence \msa{of the estimator}{} on the chosen variables for readability.
\msa{Plugging in the relations $z_\calT = \gamma \Pi_0(u_\calT - L T_\calT - \Phi_{0,\calT})^+$ and $s_\calT = k \Pi_0 T_\calT$ yields the eventually computable}{} 
\begin{eq*}
\eta_{\gamma,\calT}^2 &= \frac{1}{2}\|\nabla u_\calT - p_\calT\|_{L^2(\Omega;\R^d)}^2 + \frac{1}{2}\|\nabla T_\calT - q_\calT\|_{L^2(\Omega;\R^d)}^2 + \frac{k}{2}\|T_\calT - \Pi_0 T_\calT \|_{L^2(\Omega)}^2\\
&+ \frac{\gamma}{2}\|(u_\calT - \Phi_{0,\calT} - LT_\calT)^+ - \Pi_0 (u_\calT - \Phi_{0,\calT} - LT_\calT)^+\|_{L^2(\Omega)}^2\\
&+ \gamma (\Pi_0 (u_\calT - \Phi_{0,\calT} - LT_\calT)^+,(\Phi_{0,\calT} + LT_\calT - u_\calT)^+)_{L^2(\Omega)}.
\end{eq*}
In the remainder of this subsection, we apply our combined findings and Algorithm \ref{alg:afem}. For the calculation of the derivative, we use the expression derived in \refer{lem:apx:derivative}, but instead of the optimal solution, which we avoided to calculate, we simply plug in our choice instead and obtain
\begin{eq*}
\frac{d}{d\gamma}\eta_{\gamma,\calT}^2 \approx \frac{1}{2}\left\|(u_\calT - \Phi_{0,\calT} - LT_\calT)^+ - \Pi_0 (u_\calT - \Phi_{0,\calT} - LT_\calT)^+ \right\|_{L^2(\Omega)}^2.
\end{eq*}
Since our choice is derived from the first order system this \msa{is}{} a reasonable approximation and in fact, our expression is guaranteed to be non-negative.\\
In the following, we depict the numerical results for the thermoforming example. For the experiments we use the Doerfler marking technique with parameter $\theta = 0.1$ together with the parameters $c_\eta = 0.5$ and $c_\gamma = 0.25$. \msa{In \refer{fig:qvi_thermo:afem}}{} the behavior is tracked with respect to the degrees of freedom and the regularization parameter $\gamma$ and compared to the uniform case. The mesh as well as the comparison of the final mould $\Phi = \Phi_0 + LT$ and the membrane are depicted in \refer{fig:qvi_thermo:mesh_plot}.
\msa{Therein, the influence of the update of the regularization parameter can be clearly seen in an increase of the associated estimator. It is clearly visible in the plot, that the convergence rate of the uniform refinement follows the rate $\frac{1}{2}$, whereas the one for the adaptive case seems to follow the rate $\frac{2}{3}$. This is counter-intuitive, as first order methods usually converge maximally with rate $\frac{1}{2}$ as predicted by the a priori estimates using interpolation operators (cf. \cite[Theorem 4.4.4]{bib:BrennerScott}). This rate is already achieved by the uniform refinement strategy, and an additional increase is therefore not expected. However, it is in principle not forbidden to achieve a higher rate. An explanation might be found in the active sets and the  mesh refinement. The active sets are not depicted explicitely, but the final mould and the initial mold are depicted in \ref{fig:qvi_thermo:mesh_plot} (right). The active set can well be recognized by the scarcely refined areas in the mesh plot (right). These coincide mainly with areas, where the initial mold is affine and thus can be recovered exactly using functions in $S^1(\calT)$. This might explain the additional convergence speed.}{}

\subsection{A QVI with Application to Biomedicine}\label{ssec:qvi_membrane}
As a second example we propose the following quasi-variational inequality with relation to \msa{a}{} biomedical \msa{application}{}. The pericardium of the human heart consists of two thin layers of tissue containing a liquid between them. On both of them, forces are acting and \msa{their}{} position restricts the possible configurations of the other tissue layer. Moreover, none of the liquid is allowed to leave or enter the area between the tissues. To investigate this application, we consider the following simplified model:\\
Two membranes on an open, bounded domain $\Omega \subseteq \R^2$ are fixed on its boundary at the value zero \msa{and}{} are represented by two functions $u_1,u_2 \in H_0^1(\Omega)$, where $u_1$ models the lower membrane and $u_2$ the upper one. Of course, these two membranes are not able to perforate each other and hence the constraint $u_2 \geq u_1$ must be fulfilled almost everywhere on $\Omega$. Moreover, to model the contained liquid we require the \msa{volume}{} enclosed by the membranes to remain a predefined constant.  Without loss of generality the condition $\int_\Omega (u_2 - u_1 - 2)\dx = 0$ is proposed.
Furthermore, the forces acting on these membranes $f_1,f_2$ might depend on the positions of these membranes as well. In close proximity to the usual obstacle problem we propose the following QVI:\\
Seek $u_1,u_2 \in H_0^1(\Omega)$ that fulfil the following coupled system of variational inequalities
\begin{sys}[\label{sys:qvi_membrane}]
f_1(u) \in -\Delta u_1 + N_{K_1(u_2)}(u_1) \text{ in } H^{-1}(\Omega)\label{sys:qvi_membrane:one},\\
f_2(u) \in -\Delta u_2 + N_{K_2(u_1)}(u_2) \text{ in } H^{-1}(\Omega)\label{sys:qvi_membrane:two},
\end{sys}
with the force terms $f_1,f_2: H_0^1(\Omega)^2 \to L^2(\Omega)$. The constraint sets $K_i(u_{-i})$ read as
\begin{eq*}
\scalebox{0.96}{%
$K_i(u_{-i}) := \left\{\vphantom{\int_\Omega} u_i \in H_0^1(\Omega) : u_2 - u_1 \geq 0 \text{ a.e. on } \Omega  \text{ and } \int_\Omega (u_2 - u_1) \dx = 2|\Omega| \right\}$}
\end{eq*}
for $i = 1,2$ \msa{with}{} $u_{-1} := u_2, u_{-2} := u_1$.
The QVI \refer{sys:qvi_membrane} can be interpreted as the first order system
of the equilibrium of
\begin{eq*}
\calE_{\mathrm{QVI}}(v_1,v_2;u_1,u_2) :&= \frac{1}{2}\|\nabla v_1\|_{L^2(\Omega;\R^d)}^2 - (f_1(u),v_1)_{L^2(\Omega)} + I_{K_1(u_{\msa{2}})}(v_1)\\
&+ \frac{1}{2}\|\nabla v_2\|_{L^2(\Omega;\R^d)}^2 - (f_2(u),v_2)_{L^2(\Omega)} + I_{K_2(u_{\msa{1}})}(v_2).
\end{eq*}
In fact, one can characterize the set-valued mappings $K_i$ for $i = 1, 2$ via the \emph{joint constraint set}
\begin{eq*}
\scalebox{0.94}{%
$\calF:= \left\{ (u_1,u_2) \in H_0^1(\Omega; \R^2) : u_2 \geq u_1 \text{ a.e. on } \Omega \text{ with } \int_\Omega (u_2 - u_1 - 2) \dx = 0\right\},$}
\end{eq*}
which allows the relation $u_j \in K_j(u_{-j})$, if and only if $(u_1,u_2) \in \calF$. Based on this observation one can \msa{---}{} in close proximity to the concept of normalized/variational equilibria (see \cite{bib:Rosen} and \cite[Theorem 5]{bib:KanzowFacchinei}) \msa{---}{} formulate the equilibrium functional
\begin{eq*}
\calE_{\msa{\calF}}(v_1,v_2; u_1,u_2) &:= \frac{1}{2}\|\nabla v_1\|_{L^2(\Omega;\R^d)}^2 - (f_1(u),v_1)_{L^2(\Omega)}\\
&+ \frac{1}{2}\|\nabla v_2\|_{L^2(\Omega;\R^d)}^2 - (f_2(u),v_2)_{L^2(\Omega)} + I_{\calF}(v)
\end{eq*}
instead. With the \msa{arguments}{} therein it is straightforward to prove, that every equilibrium of $\calE$ is also an equilibrium of $\calE_{\mathrm{QVI}}$.\\
It is useful to reformulate the functional in the following way:
For a pair of membranes $u_1,u_2 \in H_0^1(\Omega)$ we define the middle position $m \in H_0^1(\Omega)$ and the halved difference $\delta \in H_0^1(\Omega)$ \msa{by}{}
\begin{eq*}
m := \frac{1}{2}(u_1 + u_2) \text{ and } \delta := \frac{1}{2}(u_2 - u_1).
\end{eq*}
Thus, the reconstruction reads $u_1 = m - \delta$ and $u_2  = m + \delta$. Setting for \msa{arbitrary}{} $v \in \calF$ the objects $m' = \frac{1}{2}(v_1 + v_2)$ and $\delta' := \frac{1}{2}(v_2 - v_1)$ we rewrite
\begin{eq*}
\calE_{\msa{\calF}}(v_1,v_2;u_1,u_2) &= \frac{1}{2}\|\nabla v_1\|_{L^2(\Omega;\R^d)}^2 + \frac{1}{2}\|\nabla v_2\|_{L^2(\Omega;\R^d)}^2 - (f_1(u),v_1)_{L^2(\Omega)}\\
&- (f_2(u),v_2)_{L^2(\Omega)} + I_\calF(v)\\
&= \frac{1}{2}\|\nabla m' - \nabla \delta'\|_{L^2(\Omega;\R^d)}^2 + \frac{1}{2}\|\nabla m' + \nabla \delta'\|_{L^2(\Omega;\R^d)}^2\\
&- (f_1(u),m' - \delta')_{L^2(\Omega)} - (f_2(u),m' + \delta')_{L^2(\Omega)}\\
&+ I_\calF(m'-\delta',m'+\delta')\\
&= \|\nabla m'\|_{L^2(\Omega;\R^d)}^2 + \|\nabla \delta'\|_{L^2(\Omega;\R^d)}^2\\
&- (f_1(m-\delta,m+\delta) + f_2(m-\delta,m+\delta),m')_{L^2(\Omega)}\\
&- (f_2(m-\delta,m+\delta) - f_1(m-\delta,m+\delta),\delta')_{L^2(\Omega)}\\
&+ I_\calF(m'-\delta',m'+\delta').
\end{eq*}
Introducing the force terms $f_m, f_\delta: H_0^1(\Omega;\R^2) \rightarrow L^2(\Omega)$
defined by
\begin{eq*}
f_m(m,d) &= \frac{1}{2}(f_1(m-\delta,m+\delta) + f_2(m-\delta,m+\delta)) \text{ and}\\
f_\delta(m,\delta) &= \frac{1}{2}(f_2(m-\delta,m+\delta)-f_1(m-\delta,m+\delta))
\end{eq*}
as well as the set
\begin{eq*}
\calG := \left\{\delta \in H_0^1(\Omega) : \delta \geq 0 \text{ a.e. on } \Omega \text{ and } \int_\Omega (\delta - 1) \dx = 0\right\}
\end{eq*}
we introduce \msa{after division by two}{} the functional
\begin{eq*}
\calE(m',\delta';m,\delta) :&= \frac{1}{2}\|\nabla m'\|_{L^2(\Omega;\R^d)}^2 - (f_m(m,\delta),m')_{L^2(\Omega)}\\
&+ \frac{1}{2}\|\nabla \delta'\|_{L^2(\Omega;\R^d)}^2 - (f_\delta(m,\delta),\delta')_{L^2(\Omega)} + I_\calG(\delta').
\end{eq*}
In the following, we seek to work with this equilibrium formulation for several reasons\msa{:}{}\\
On the one hand, we eliminated the dependence of the indicator functional on the feedback component. Hence, we are enabled to utilize the slightly more convenient framework of \refer{thm:myreg:mosco} and obtain even Mosco-convergence results. On the other hand, \msa{we notice, that}{} this approach is in principle not limited to this \msa{particular}{} formulation: One \msa{notes}{}, that the minimization problem induced by $\calE_{\mathrm{QVI}}$ decomposes \msa{into}{} two minimization problems for $v_1,v_2$ respectively. Using \msa{an interpolation}{} with \msa{respect to}{} a parameter $\lambda \in (0,1)$ one \msa{derives}{} an analogous \msa{functional}. \msa{The latter is used to}{} prove, that in many cases the QVI in \refer{sys:qvi_membrane} has infinitely many solutions.
\begin{thm}\label{thm:qvi_membrane:infinite_solutions}
\msa{Assume, that}{} the functional $f_\delta$ only depends on $\delta$. \msa{Consider the following VI}:\\
Seek $\delta \in H_0^1(\Omega)$ such that
\begin{eq}[\label{eq:qvi_membrane:delta_reduced}]
0 \in -\Delta \delta - f_\delta(\delta) + N_\calG(\delta) \text{ in } H^{-1}(\Omega),
\end{eq}
as well as the \msa{following}{} PDE\msa{:\\
For given $\xi \in H^{-1}(\Omega)$ and $\delta \in \calG$ seek $m \in H_0^1(\Omega)$ such that
\begin{eq}[\label{eq:qvi_membrane:m_reduced}]
-\Delta m - f_m(m,\delta) &= \xi &&\text{in } \Omega,\\
m &= 0 &&\text{on } \partial \Omega,
\end{eq}
is fulfilled.}{} \msa{If  \refer{eq:qvi_membrane:delta_reduced} has a solution  and \refer{eq:qvi_membrane:m_reduced} admits a solution for every given $\delta \in \calG$ and $\xi \in H^{-1}(\Omega)$, then}{} the \msa{QVI}{} \refer{sys:qvi_membrane} has a solution.\\
Moreover, if there exists a solution $\delta \in H_0^1(\Omega)$ of \refer{eq:qvi_membrane:delta_reduced} with  $-\Delta \delta - f_\delta(\delta) \neq 0$, then the \msa{QVI}{} \refer{sys:qvi_membrane} has infinitely many solutions.
\end{thm}
\begin{proof}
Take an arbitrary $\lambda \in (0,1)$. Considering the \msa{QVI}{} \refer{sys:qvi_membrane} define the functional
\begin{eq*}
\calE_{QVI}^\lambda(v_1,v_2;u_1,u_2) := \lambda\left(\frac{1}{2}\|\nabla v_1\|_{L^2(\Omega;\R^d)}^2 - (f_1(u),v_1)_{L^2(\Omega)} + I_{K_1(u_2)}(v_1)\right)\\
+ (1-\lambda)\left(\frac{1}{2}\|\nabla v_2\|_{L^2(\Omega;\R^d)}^2 - (f_2(u),v_2)_{L^2(\Omega)} + I_{K_2(u_1)}(v_2)\right).
\end{eq*}
Here, the part of the functional minimized with respect to $v_1$ has been scaled with $\lambda$ and the part with respect to $v_2$ has been scaled with $(1-\lambda)$. Again, every equilibrium of $\calE_{QVI}^\lambda$ is a solution of \refer{sys:qvi_membrane}. Hence, we proceed as before with the derivation of the corresponding variational equilibrium by substituting the indicators by $I_\calF(v)$ and define
\begin{eq*}
&\calE^\lambda(\delta',m';\delta,m) := \lambda\left(\frac{1}{2}\|\nabla m' - \nabla \delta\|_{L^2(\Omega;\R^d)}^2 - (f_m(m,\delta) + f_\delta(m,\delta),m' - \delta')_{L^2(\Omega)}\right)\\
&+ (1-\lambda)\left(\frac{1}{2}\|\nabla m' + \nabla \delta'\|_{L^2(\Omega;\R^d)}^2 - (f_m(m,\delta) + f_\delta(m,\delta),m'+\delta')_{L^2(\Omega)}\right) + I_\calG(\delta')\\
&= \frac{1}{2}\|\nabla m'\|_{L^2(\Omega;\R^d)}^2 + \frac{1}{2}\|\nabla \delta'\|_{L^2(\Omega;\R^d)} + (1-2\lambda)(\nabla m', \nabla \delta')_{L^2(\Omega;\R^d)}\\
&- (f_m(m,\delta), m')_{L^2} + (2\lambda - 1)(f_m(m,\delta), \delta')_{L^2}\\
&- (f_\delta(m,\delta),\delta')_{L^2} + (2\lambda - 1)(f_\delta(m,\delta),m')_{L^2} + I_\calG(\delta').
\end{eq*}
\msa{Then, the}{} first order system \msa{of $\calE^\lambda$}{} reads as
\begin{eq}[\label{eq:qvi_thermo:lambda:firstorder}]
0 &= -\Delta m - (1-2\lambda)\Delta \delta - f_m(m,\delta) - (1-2\lambda)f_\delta(\delta)\\
&= -\Delta m - f_m(m,\delta) + (1-2\lambda)( -\Delta \delta - f_\delta(\delta)) \text{ in } \Omega, 
m = 0 \text{ on } \partial \Omega.\\
0&\in -\Delta \delta - (1-2\lambda)\Delta m - f_\delta(\delta) - (1-2\lambda)f_m(m,\delta) + N_\calG(\delta)\\
&= -\Delta \delta - f_\delta(\delta) + (1-2\lambda)(-\Delta m - f_m(m,\delta)) + N_\calG(\delta) 
\end{eq}
Plugging in the first equation into the VI leads to 
\begin{eq*}
0 &\in -\Delta \delta - f_\delta(\delta) - (1-2\lambda)^2(-\Delta \delta - f_\delta(\delta)) + N_\calG(\delta)\\
&= 4\lambda(1-\lambda)(-\Delta \delta - f_\delta(\delta)) + N_\calG(\delta).
\end{eq*}
Dividing by $4\lambda(1-\lambda) > 0$ leads to \refer{eq:qvi_membrane:delta_reduced}. \msa{By assumption, this VI has a solution and solving \refer{eq:qvi_membrane:m_reduced} for $\xi = -(1-2\lambda)(-\Delta \delta - f_\delta(\delta))$ yields the existence of a solution of \refer{eq:qvi_thermo:lambda:firstorder}, and thus $\calE^\lambda$ has an equilibrium.}{}\\
If moreover, $-\Delta \delta - f_\delta(\delta) \neq 0$ holds, then the right hand side in \msa{\refer{eq:qvi_membrane:m_reduced} differs with changing $\lambda$ for the previous choice of $\xi$.}{}
\msa{As}{} every solution mapping $\lambda \mapsto m_\lambda$ is injective \msa{and every}{} pair $(\delta,m_\lambda)$ solves the \msa{QVI}{} \refer{sys:qvi_membrane}, \msa{there}{} exist infinitely many solutions. 
\end{proof}

\begin{figure}
\begin{tabular}{p{0.49\textwidth}p{0.49\textwidth}}
\begin{tikzpicture}[scale = 4.0]
	\draw[->, color = gray] (-0.1,0) -- (1.2,0);
	\draw[->, color = gray] (0,-0.1) -- (0,1.2);
	\draw (1.0,-0.025) -- (1.0,0.025);
	\draw (-0.025,1.0) -- (0.025,1.0);
	\draw (0.5,-0.025) -- (0.5,0.025);
	\draw (-0.025,0.5) -- (0.025,0.5);
	\draw (1/6,-0.025) -- (1/6,0.025);
	\draw (1/3,-0.025) -- (1/3,0.025);
	\draw (-0.025,2/3) -- (0.025,2/3);
	\draw (-0.025,5/6) -- (0.025,5/6);
	\node[anchor = east] at (-0.025,2/3) {$\frac{2}{3}$};
	\node[anchor = east] at (-0.025,5/6) {$\frac{5}{6}$};
	\node[anchor = east] at (-0.025,0.5) {$\frac{1}{2}$};
	\node[anchor = north] at (0.5,-0.025) {$\frac{1}{2}$};
	\node[anchor = north] at (1.0,-0.025) {1};
	\node[anchor = north east] at (-0.025,-0.025) {0};
	\node[anchor = east] at (-0.025,1.0) {1};
	\node[anchor = north] at (1/6,-0.025) {$\frac{1}{6}$};
	\node[anchor = north] at (1/3,-0.025) {$\frac{1}{3}$};
	\draw (0,0) -- (0,1.0) -- (1.0,1.0) -- (1.0,0.5) -- (0.5,0.5) -- (0.5,0.0) -- cycle;
	\draw[fill, black!20] (1/6,1/6) -- (1/3,1/6) -- (1/3,2/3) -- (5/6,2/3) -- (5/6,5/6) -- (1/6,5/6) -- cycle;
\end{tikzpicture}
&
\begin{tikzpicture}[scale = 4.0]
	\draw[->, color = gray] (-0.1,0) -- (1.2,0);
	\draw[->, color = gray] (0,-0.1) -- (0,1.2);
	\draw[fill, black!20] (0,0) -- (0,1.0) -- (1.0,1.0) -- (1.0,0) -- cycle;
	\draw[fill, white] (1/4,1/4) -- (1/4,3/4) -- (3/4,3/4) -- (3/4,1/4) -- cycle;
	
	\draw (1.0,-0.05) -- (1.0,0.05);
	\draw (-0.05,1.0) -- (0.05,1.0);
	\draw (0.5,-0.025) -- (0.5,0.025);
	\draw (0.25,-0.025) -- (0.25,0.025);
	\draw (0.75,-0.025) -- (0.75,0.025);
	\draw (-0.025,0.5) -- (0.025,0.5);
	\draw (-0.025,0.25) -- (0.025,0.25);
	\draw (-0.025,0.75) -- (0.025,0.75);
	\node[anchor = north] at (1.0,-0.05) {1};
	\node[anchor = north east] at (-0.05,-0.05) {0};
	\node[anchor = east] at (-0.05,1.0) {1};
	\node[anchor = north] at (0.5,-0.025) {$\frac{1}{2}$};
	\node[anchor = north] at (0.25,-0.025) {$\frac{1}{4}$};
	\node[anchor = north] at (0.75,-0.025) {$\frac{3}{4}$};
	\node[anchor = east] at (-0.025,0.5) {$\frac{1}{2}$};
	\node[anchor = east] at (-0.025,0.25) {$\frac{1}{4}$};
	\node[anchor = east] at (-0.025,0.75) {$\frac{3}{4}$};
	\draw (0,0) -- (0,1.0) -- (1.0,1.0) -- (1.0,0.5) -- (0.5,0.5) -- (1.0,0.5) -- (1.0,0.0) -- cycle;
\end{tikzpicture}
\end{tabular}
\caption{Depiction of the domain $\Omega$, Left: L-shape. Right: Slit. In gray: area, where the forces $f_{\mathrm{const},1}, f_{\mathrm{const},2}$ take action.}\label{fig:qvi_membrane:grayareas}
\end{figure}
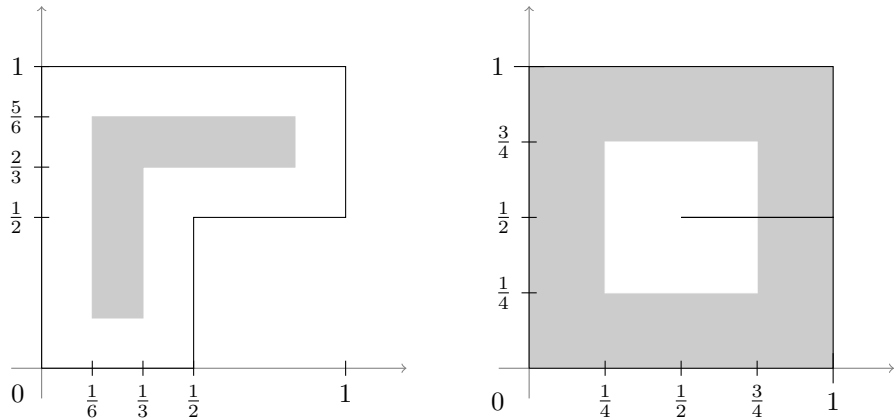

\begin{figure}[h]
\includegraphics[trim = {890, 250, 890, 250}, clip, height = 0.175\textheight]{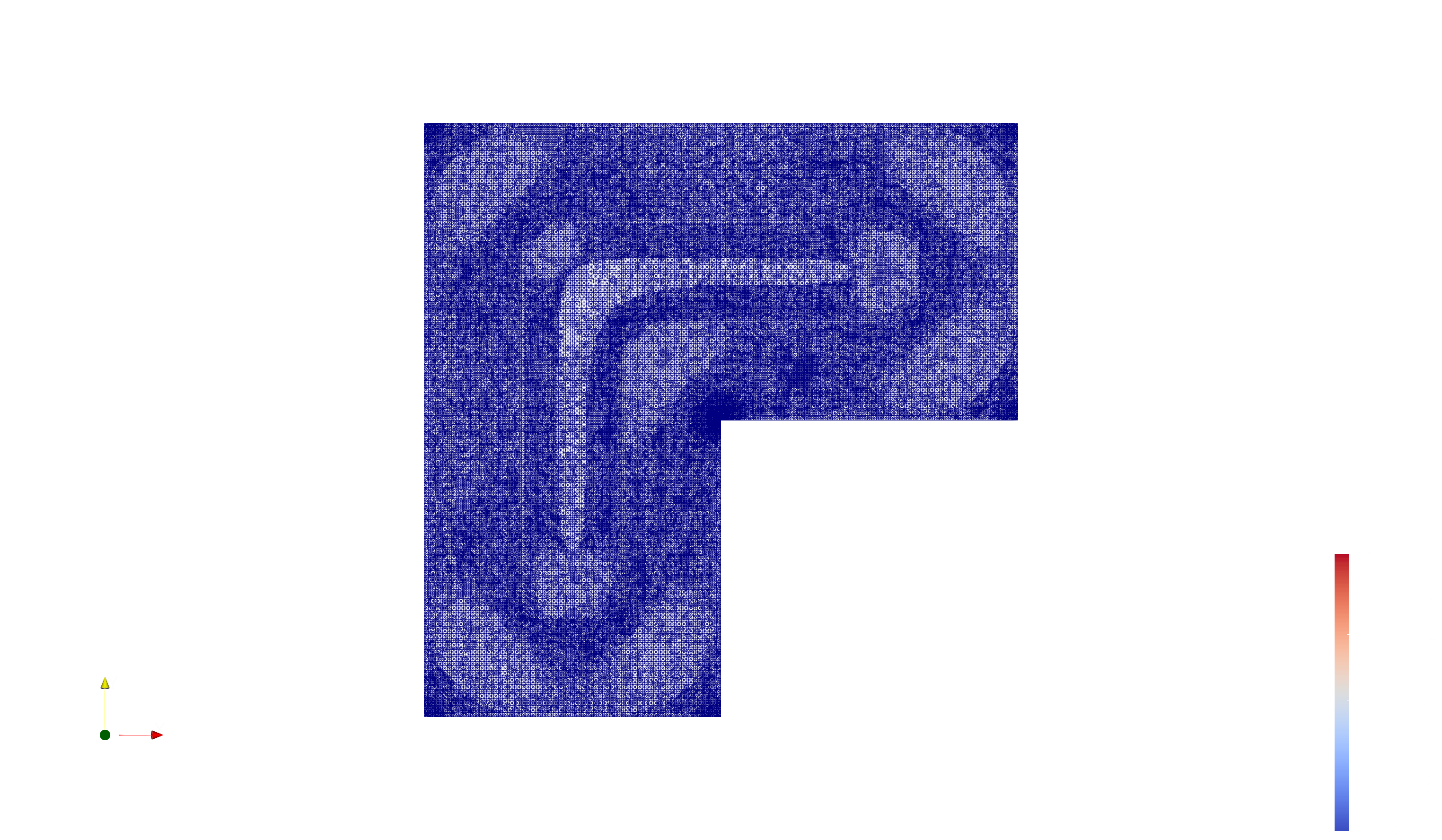}
\includegraphics[trim = {650, 500, 300, 420}, clip, height = 0.175\textheight]{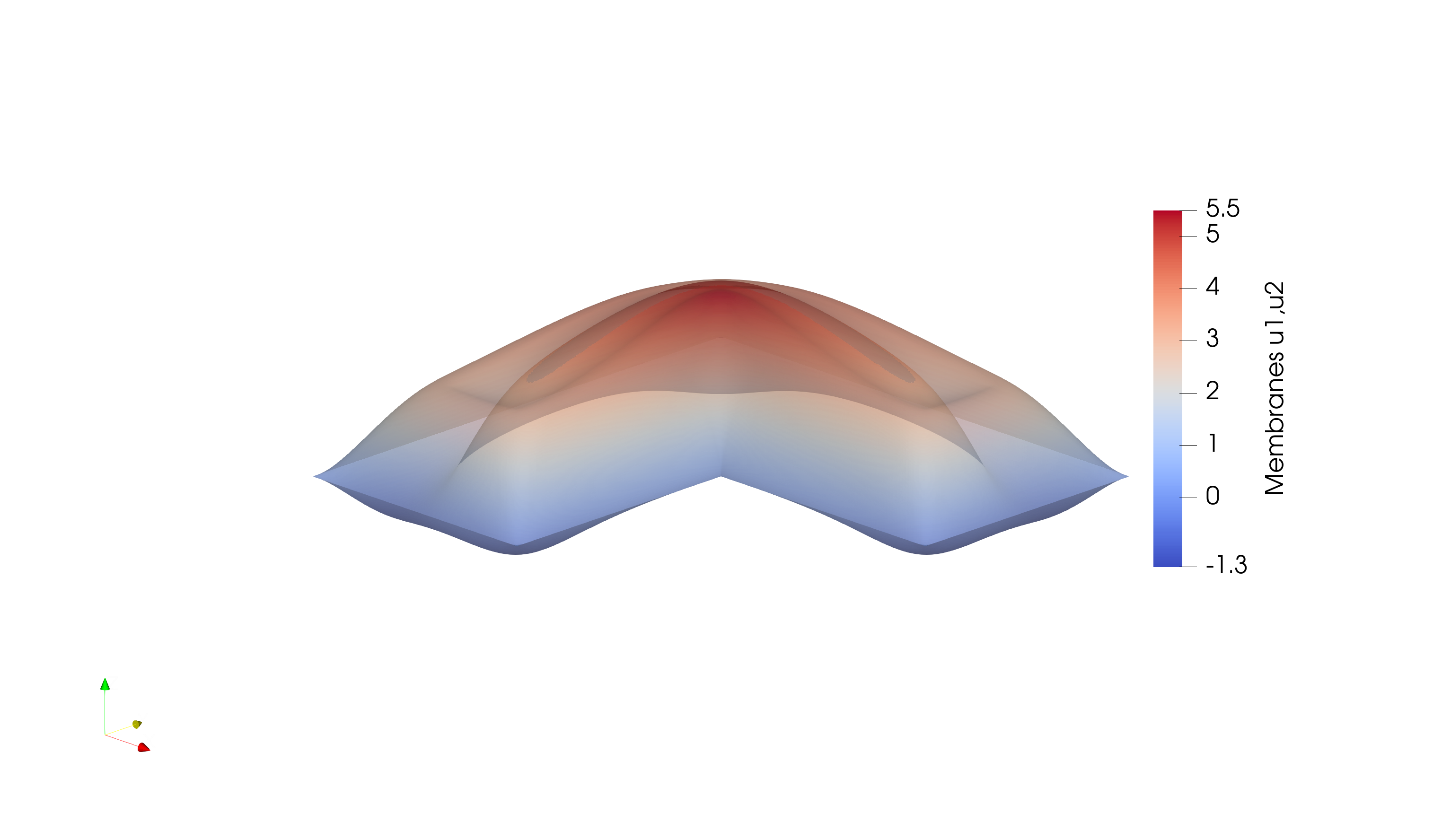}
\caption{Left: Adaptively refined mesh. Right: Plot of the two membranes $u_1$ (lower) and $u_2$ (upper).}\label{fig:qvi_membrane:Lshape:mesh_plot}
\end{figure}

\begin{figure}[h]
\includegraphics[trim = {890, 250, 890, 250}, clip, height = 0.175\textheight]{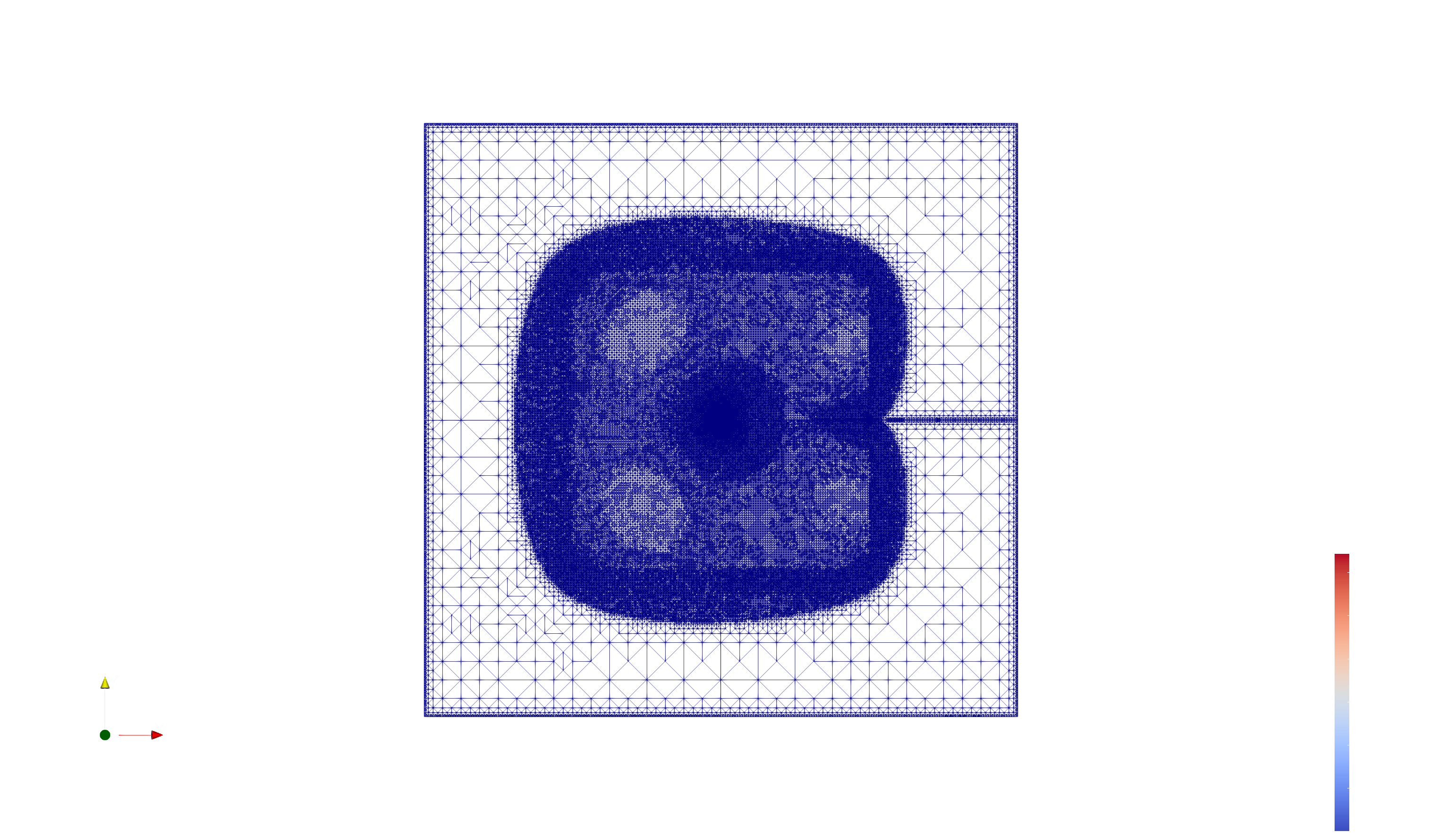}
\includegraphics[trim = {650, 500, 300, 450}, clip, height = 0.175\textheight]{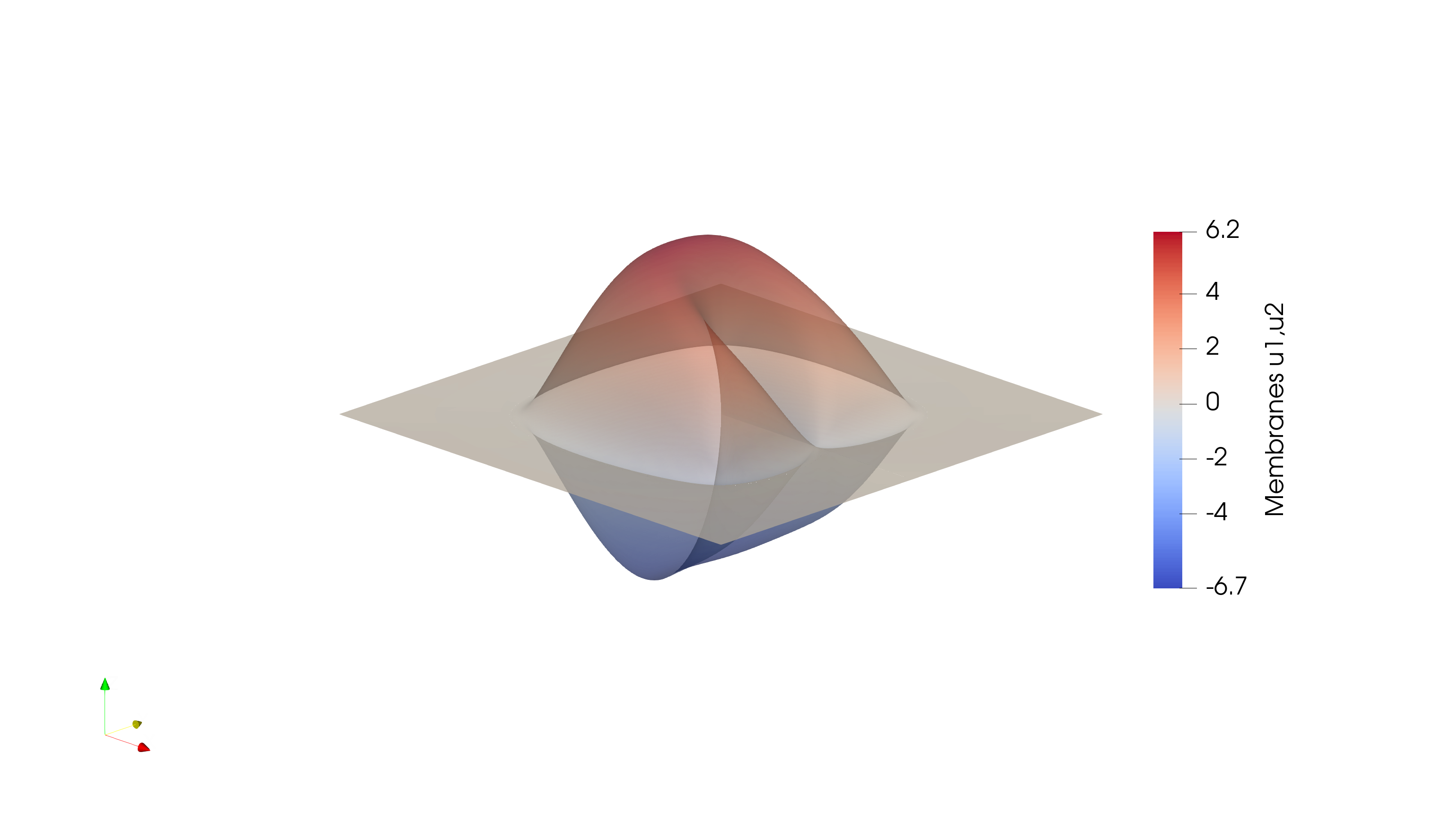}
\caption{Left: Adaptively refined mesh. Right: Plot of the two membranes $u_1$ (lower) and $u_2$ (upper).}\label{fig:qvi_membrane:slit:mesh_plot}
\end{figure}

On the one hand, \refer{thm:qvi_membrane:infinite_solutions} proves, that in this particular setting the existence question of equilibria of $\calE$ boils down to the existence of a solution of the VI in \refer{eq:qvi_membrane:delta_reduced}.  The above conditions are especially fulfilled for the simplest possible case $f_1 \equiv f_2 \equiv 0$. \msa{Hence,}{} one cannot expect uniqueness of solutions for the QVI in general. On the other hand\msa{,}{} it might be, that \msa{$\calE$}{} admits a unique solution nevertheless.\\
The first order system of equilibria of $\calE$ reads as the system
\begin{sys*}
-\Delta m - f_m(m,\delta) &= 0 \text{ in } \Omega,\\
m &= 0 \text{ on } \partial \Omega,\\
f_\delta(m,\delta) &\in -\Delta \delta + N_\calG(\delta) \text{ in } H^{-1}(\Omega),
\end{sys*}
which consists of a partial differential equation coupled with a variational inequality. We propose the penalized equilibrium problem induced by the functional
\begin{eq*}
\calE_{\msa{\gamma}}(m',&\delta';m,\delta) := \frac{1}{2}\|\nabla m'\|_{L^2(\Omega;\R^d)}^2 - (f_m(m,\delta),m')_{L^2(\Omega)}\\
&+ \frac{1}{2}\|\nabla \delta'\|_{L^2(\Omega;\R^d)}^2 - (f_\delta(m,\delta),\delta')_{L^2(\Omega)} + \frac{\gamma}{2}\int_\Omega (-\delta')^{2+} \dx + I_V(\delta'),
\end{eq*}
where $V$ is the affine subset 
\begin{eq*}
V:= \left\{\delta \in H_0^1(\Omega): \int_\Omega (\delta - 1) \dx = 0\right\}.
\end{eq*}
Here, we used $\calG = L^2_+(\Omega) \cap V$, respectively the decomposition $I_\calG = I_{L^2_+(\Omega)} + I_V$. Defining the set $V_{L^2} := \{\delta \in L^2(\Omega) \msa{:}{} \int_\Omega (\delta - 1) \dx = 0\}$ we rewrite $I_V(\ccdot) = I_{V_{L^2}}(i_{L^2} (\ccdot))$. As it can be seen in the formulation of the functional we did indeed use a Moreau--Yosida regularization for the sign condition on $\delta$, but we left the volume constraint untouched. There are several reasons to do so: On the one hand, the application of a Moreau--Yosida on $V_{L^2}$ leads to the term $\frac{\gamma}{2}\left( \int_\Omega (\delta - 1) \dx \right)^2$. This might lead to difficulties, when it is tried to formulate an error estimator, since it is not easily decomposed into a sum over triangles $T \in …\calT$. Moreover, it is not clear to which extent a penalty for \msa{a}{} global constraint can provide local information. On the other hand, there is no necessity to regularize the volume constraint.  In fact, the proposition of the obstacle-type constraint does induce difficulties in the fulfilment of the constraint qualification, whereas the volume constraint does not, as we see next for the derivation of the first order condition:\\
Using the developed methodology we define the spaces $U = H_0^1(\Omega) \times H_0^1(\Omega)$ and $Y := L^2(\Omega;\R^d) \times L^2(\Omega) \times L^2(\Omega;\R^d) \times L^2(\Omega) \times L^2(\Omega)$ as well as $Z:= L^2(\Omega)$. As functionals we define $\msa{F(m,\delta)}{} = 0$ and
\begin{eq*}
G(p_m,p_\delta,z_m,&z_\delta,z_V;m,\delta) := \frac{1}{2}\|p_m\|_{L^2(\Omega;\R^d)}^2 + \frac{1}{2}\|p_\delta\|_{L^2(\Omega;\R^d)}^{\msa{2}}\\
&- (f_m(m,\delta),z_m)_{L^2(\Omega)} - (f_\delta(m,\delta),z_\delta)_{L^2(\Omega)} + I_{V_{L^2}}(z_V)
\end{eq*}
as well as $H(z) := I_{L^2_+(\Omega)}(z)$. As linear operators we choose $A(m,\delta) := (\nabla m, \nabla \delta, i_{L^2(\Omega)}m, i_{L^2(\Omega)}\delta,i_{L^2(\Omega)}\delta)$ and $B(m,\delta)  := i_{L^2(\Omega)} \delta$. The constraint qualification \refer{eq:duality:my:cq} reads equivalently as
\begin{eq*}
0 \in \core{AU - L^2(\Omega;\R^d) \times L^2(\Omega;\R^d) \times L^2(\Omega) \times L^2(\Omega) \times V}.
\end{eq*}
It is sufficient to prove, that the set inside the core is indeed the whole space $Y$. Therefore, let arbitrary $q_m, q_\delta \in L^2(\Omega;\R^d)$ and $y_m, y_\delta,y_V \in L^2(\Omega)$ be given and we seek to find $m' \in H_0^1(\Omega)$ and $\delta' \in V$ together with $p_m,p_\delta \in L^2(\Omega;\R^d)$ and $z_m,z_\delta,z_V \in L^2(\Omega)$\msa{, such that the following system is fulfilled
\begin{eq*}
\nabla m - p_m &= q_m,\\
\nabla \delta - p_\delta &= q_\delta,\\
i_{L^2(\Omega)}m - z_m &= y_m,\\
i_{L^2(\Omega)}\delta - z_\delta &= y_\delta,\\
i_{L^2(\Omega)}\delta - z_V &= y_V.
\end{eq*}}
\msa{For that, take}{} $\phi \in H_0^1(\Omega)$ with $\int_\Omega \phi\, \dx = 1$ and for \msa{given}{} $\msa{y_V}{} \in L^2(\Omega)$ define $\delta' := y_V - (\int_\Omega (y_V -1) \dx)\cdot \phi$ and see
\begin{eq*}
\int_\Omega (\delta' - 1) \dx = \int_\Omega (y_V - 1) \dx - \int_\Omega (y_V - 1) \dx \cdot \int_\Omega \phi \dx = 0
\end{eq*}
\msa{or in other words $\delta' \in V_{L^2}$. Then, we}{} obtain the decomposition $y_V = \delta' - z_V$ with $z_V = -\int_\Omega (y_V - 1) \dx \cdot \phi$. Thus, we can furthermore take $z_\delta = y_\delta - \delta'$ as well as $m' = 0$ and $p_m = q_m$ and $p_\delta = q_\delta$ and obtain the assertion.\\
After the discussion of the constraint qualification we \msa{apply}{} the sum rule \msa{to derive}{} the first order \msa{condition}{}. Therefore, \msa{we show for}{} the normal cone $N_{V_{L^2}}(v) = \R$ for all $v \in V_{L^2}$. To see this, let \msa{first}{} $\mu \in L^2(\Omega)$ be an element in the normal cone of $V_{L^2}$. \msa{By its definition holds}{}
\begin{eq*}
(\mu, v' - v)_{L^2(\Omega)} \leq 0 \text{ for all } \msa{v}' \in V_{L^2}.
\end{eq*}
Taking $v' := v + \mu - \frac{1}{\operatorname{area}(\Omega)}\int (v + \mu - 1)\dx = v + \mu - \frac{1}{\operatorname{area}(\Omega)}\int_\Omega \mu \dx \in V_{L^2}$ we obtain
\begin{eq*}
0 \geq (\mu, v' - v)_{L^2(\Omega)} = \left(\mu, \mu - \fint_\Omega \mu \dx\right) = \left\|\mu - \fint_\Omega \mu \dx\right\|_{L^2(\Omega)}^2 \geq 0
\end{eq*}
and hence $\mu  =\fint_\Omega \mu \dx$, which implies $\mu \in \R$.\\
\msa{For the other direction, observe, that}{} for all $\mu \in \R$ holds $(\mu, v' - v)_{L^2(\Omega)} = \mu (\int_\Omega (v' - 1) \dx - \int_\Omega (v - 1) \dx ) = 0$ for all $v,v' \in V_{L^2}$ and hence $N_{V_{L^2}} \equiv \R$.\\
Thus, the equilibrium of the penalized problem fulfils the following system: Seek $m_\gamma,\delta_\gamma \in H_0^1(\Omega)$ and $\mu_\gamma \in \R$ such that
\begin{sys}[\label{sys:qvi_membrane:firstorder}]
-\Delta m_\gamma &= f_m(m_\gamma,\delta_\gamma) &\text{in } \Omega,&& m_\gamma &= 0 \text{ on } \partial \Omega,\\
-\Delta \delta_\gamma - \gamma(-\delta_\gamma)^+ + \mu_\gamma &= f_\delta(m_\gamma,\delta_\gamma) &\text{in } \Omega,&& \delta_\gamma &= 0 \text{ on } \partial \Omega,\\
\int_\Omega (\delta_\gamma - 1) \dx &= 0
\end{sys}
hold.\\
\msa{To extract weakly convergent subsequences of equilibria as requested in \refer{prop:equi:gamma:mini_convergence} and \refer{thm:diag} we assure the boundedness of sequences of equilibria. For the derivation of a priori bounds, we propose the following assumptions:}{}\\
The operators $f_\delta,f_m$ only depend on $\delta \in H_0^1(\Omega)$\msa{, the operator $f_m$ is bounded, i.e. mapping bounded subsets to bounded subsets,}{} and the operator $f_\delta$ can be decomposed as $f_\delta = g_\delta + h_\delta$\msa{,}{} such that $-\Delta - g_\delta$ is strongly monotone and the operator $h_\delta: H_0^1(\Omega) \rightarrow L^2(\Omega)$ is uniformly bounded on $\calG$ by a constant $B_h > 0$. \msa{Then}, a hierarchy is established\msa{: One solves}{} the second equation with respect to $\delta$ \msa{first,}{} and afterwards the (linear) PDE with respect to $m$. \msa{Under these assumptions}, there exists a constant $\alpha > 0$ with
\begin{eq*}
\alpha \|\nabla \delta - \nabla \delta'\|^2_{L^2\msa{(\Omega)}{}} \leq \|\nabla \delta - \nabla \delta'\|^2_{L^2\msa{(\Omega)}} - (g_\delta(\delta) - g_\delta(\delta'),\delta - \delta')_{L^2\msa{(\Omega)}}.
\end{eq*}
Hence, we obtain for an arbitrary $\delta_0 \in \calG$ the estimate
\begin{eq*}
\alpha \|\nabla \delta_\gamma &- \nabla \delta_0\|^2_{L^2\msa{(\Omega)}} \leq \|\nabla \delta_\gamma - \nabla \delta_0\|^2_{L^2\msa{(\Omega)}} - \gamma ((-\delta_\gamma)^+ - (-\delta_0)^+, \delta_\gamma - \delta_0)_{L^2\msa{(\Omega)}}\\
&- (g_\delta(\delta_\gamma) - g_\delta(\delta_0),\delta_\gamma - \delta_0)_{L^2\msa{(\Omega)}}\\
&= -(\nabla \delta_0, \nabla \delta_\gamma - \nabla \delta_0)_{L^2\msa{(\Omega)}} + (g_\delta(\delta_0),\delta_\gamma - \delta_0)_{L^2\msa{(\Omega)}}\\
&-(h(\delta_\gamma), \delta_\gamma - \delta_0)_{L^2\msa{(\Omega)}} - \mu \int_\Omega (\delta_\gamma - \delta')\dx\\
&\leq \left(\|\nabla \delta_0\|_{L^2(\Omega\msa{;\R^d})} + \|g(\delta_0)\|_{L^2(\Omega)} + \|h(\delta_\gamma)\|_{L^2(\Omega)}\right) \|\delta_\gamma - \delta_0\|_{L^2\msa{(\Omega)}}
\end{eq*}
and \msa{using Friedrich's inequality eventually}{} an a priori bound for $\delta_\gamma$ by
\begin{eq*}
\|\nabla \delta_\gamma\|_{L^2(\Omega\msa{;\R^d})} \lesssim \|\nabla \delta_0\|_{L^2(\Omega;\R^d)} + \|g(\delta_0)\|_{L^2(\Omega)} + B_h.
\end{eq*}
Subsequently, we obtain an a priori bound for $m_\gamma$ by 
$\|\nabla m_\gamma\|_{L^2\msa{(\Omega)}} \leq \|f_m(\delta_\gamma)\|_{L^2\msa{(\Omega)}}$
and hence the boundedness of the sequence of equilibria $(m_\gamma,\delta_\gamma) \subset H_0^1(\Omega;\R^2)$. Using \refer{thm:myreg:mosco} we deduce from the Mosco-convergence the existence as well as convergence to an equilibrium of $\calE$ and hence to a solution of the QVI \refer{sys:qvi_membrane}. In fact, in the light of \refer{thm:diag} we obtain the existence of a subsequence of approximate equilibria, that converge strongly to a solution of the QVI. With \refer{lem:duality:my} together with Calculation \ref{calc:conjugate} and the decomposition used for the first order condition we \msa{derive the}{} dual functional
\begin{eq*}
\calD_\gamma(p_m,p_\delta,&z_V,z; m, \delta) := \frac{1}{2}\|p_m\|_{L^2(\Omega;\R^d)}^2 + \frac{1}{2}\|p_\delta\|_{L^2(\Omega;\R^d)}^2\\
&+ I_{\{0\}}(-\ddiv p_m - f_m(\delta)) + I_{\{0\}}(-\ddiv p_\delta + z_V + z - f_\delta(\delta))\\
&+ \mathrm{area}(\Omega) z_V + I_{L^2_-(\Omega)}(z) + \frac{1}{2\gamma}\|z\|_{L^2(\Omega)}^2.
\end{eq*}
\msa{To perform numerical computations}, we \msa{propose}{} the following situation: There exist constant forces $f_1,f_2 \in L^2(\Omega)$ solely acting on the membranes $u_1$ \msa{and}{} $u_2$\msa{, respectively}. Moreover, the fluid has a weight that is acting on the lower membrane vertically with a force $-\alpha(u_2-u_1)$ and hence we have the force operators
\msa{
\begin{eq*}
f_1(u_1,u_2) = f_{\mathrm{const},1} - \alpha (u_2-u_1) \text{ and } f_2(u_1,u_2) = f_{\mathrm{const},2},
\end{eq*}
respectively
\begin{eq*}
f_m(\delta) = \frac{1}{2}(f_{\mathrm{const},1} + f_{\mathrm{const},2}) - \alpha \delta \text{ and } f_\delta(\delta) = \frac{1}{2}(f_{\mathrm{const},2} - f_{\mathrm{const},1}) + \alpha\delta.
\end{eq*}}{}%
\msa{Then,}{} $g_\delta(\delta) = \alpha \delta$ and $h_\delta(\delta) = \frac{1}{2}(f_{\mathrm{const},2} - f_{\mathrm{const},1})$. Assuming, that the constant forces can be resolved exactly on a sufficiently fine mesh, we do not need to introduce functionals $\calE_\gamma^\calT,\calD_\gamma^\calT$ and hence do not need to discuss any oscillatory terms. \msa{Similar to the previous examples we discretize the primal objects with functions in $S_0^1(\calT)$.
We choose for given $m_\calT,\delta_\calT \in S_0^1(\calT)$ the function $p_{m,\calT}, p_{\delta,\calT} \in RT_0(\Omega)$ to be the solution of the mixed formulation of the Poisson equation with right hand side $f_m(\delta_\calT)$ respectively $-\mu_\calT + \gamma(-\delta_\calT)^+ + f_\delta(\delta_\calT)$.
Using partial integration we obtain
\begin{eq*}
(f_m(\delta_\calT),m_\calT) = (p_{m,\calT},\nabla m_\calT)_{L^2(\Omega;\R^2)}
\end{eq*}
as well as
\begin{eq*}
(f_\delta(\delta_\calT),\delta_\calT) &= (p_{\delta,\calT},\nabla \delta_\calT)_{L^2(\Omega;\R^2)} - (z_{V,\calT},\delta_\calT)_{L^2(\Omega)} - (z_\calT,\delta_\calT)_{L^2(\Omega)}\\
&= (p_{\delta,\calT},\nabla \delta_\calT)_{L^2(\Omega;\R^2)} - \mathrm{area}(\Omega) z_{V,\calT} - (z_\calT,\delta_\calT)_{L^2(\Omega)}.
\end{eq*}}{}
Using $z_\calT = -\gamma \Pi_0(-\delta_\calT)^+$ the estimator reads as
\begin{eq*}
\eta_{\gamma,\calT}^2 &= \calE_\gamma(m_{\msa{\calT}},\delta_{\msa{\calT}};m_{\msa{\calT}},\delta_{\msa{\calT}}) + \calD_\gamma(p_{m,\calT},p_{\delta,\calT},z_{V,\calT},z_\calT;m_\calT,\delta_\calT)\\
&= \frac{1}{2}\|\nabla m_\calT - p_{m,\calT}\|_{L^2(\Omega;\R^{\msa{2}})}^2 + \frac{1}{2}\|\nabla \delta_\calT - p_{\delta,\calT}\|_{L^2(\Omega;\R^2)}^2\\
&+ \frac{\gamma}{2}\int_\Omega (-\delta_\calT)^{2+} \dx - (z_\calT,\delta_\calT)_{L^2(\Omega)} + I_{L^2_-(\Omega)}(z_\calT) + \frac{1}{2\gamma}\|z_\calT\|_{L^2(\Omega)}^2\\
&= \frac{1}{2}\|\nabla m_\calT - p_{m,\calT}\|_{L^2(\Omega;\R^d)}^2 + \frac{1}{2}\|\nabla \delta_\calT - p_{\delta,\calT}\|_{L^2(\Omega;\R^2)}^2\\
&+ \frac{\gamma}{2}\| (-\delta_\calT) - \Pi_0(-\delta_\calT)^+ \|_{L^2(\Omega)}^2 + \gamma (\Pi_0 (-\delta_\calT)^+, \delta_\calT^+)_{L^2(\Omega)},
\end{eq*}
where we again suppress the dependence on the chosen objects in the notation \msa{of the estimator}.\\
\msa{Next, we perform our numerical experiments. For this sake, we propose two different settings: In the first case, we consider the \emph{L-shape domain} $\Omega := (0,1)^2\backslash \left([\frac{1}{2},1)\times (0,\frac{1}{2}]\right)$ and use $f_{\mathrm{const},1} = 1000$ and $f_{\mathrm{const},2} = -500$ taking action on the gray areas depicted in \refer{fig:qvi_membrane:grayareas}. Hence, we have there in total a force pushing the two membranes upwards.\\
In the second case, we set $\Omega = (0,1)^2 \backslash \left([\frac{1}{2},1) \times \{\frac{1}{2}\} \right)$ as the \emph{slit domain} and use $f_{\mathrm{const},1} = -f_{\mathrm{const},2} = 1000$ acting on the respective gray area in \refer{fig:qvi_membrane:grayareas} and obtain forces, that clamp both membranes together. For the gravitational force $\alpha = 2$ is taken, to ensure the strong monotony. In Algorithm \ref{alg:afem} we choose the Doerfler parameter $\theta = 0.1$ and $c_\gamma = 0.1$. Our numerical findings are depicted in the Figures \ref{fig:qvi_membrane:Lshape:afem} to \ref{fig:qvi_membrane:slit:mesh_plot}.}{}

\msa{In the plots one clearly sees the membranes being clamped together on the gray areas. Due to the volume constraint the membranes divert elsewhere from each other drastically more than it would be expected without the presence of that condition. Besides the constant forces, we also have a gravitational force pushing the lower membrane downwards. Clearly, the latter force is significantly weaker than the first one. However, its influence can be seen in the plots as well: For the L-shape domain we see areas, that go below zero. This is partly to the aforementioned volume constraint, but also influenced by the force acting downwards. The influence of the latter can be seen clearer for the slit domain. There, a brief observation suggests a symmetry, but closer inspection of the scales reveals the middle area to be shifted downwards.\\
The plots of the adaptively refined meshes reveal the structure of the gray areas, as well as the contact area of the two membranes. The latter effect is far more prominent for the slit domain.\\
Addressing the convergence rate one can see an increase in both cases as it is expected for non-convex domains. The convergence rate in the adaptive case is $\frac{1}{2}$ each. For the uniform case, we obtain the rates $\frac{2}{5}$ and $\frac{1}{4}$ for the L-shape and slit domain. The effect of the update of the regularization parameter $\gamma$ is not as profound as for the thermoforming example. This is in part related to the smaller parameter $c_\gamma$. Despite this, we obtain a similar profile regarding the update of the regularization parameter. This indicates an overall less sensitive dependence of the estimators on $\gamma$.}{}

\begin{figure}[h]
\newcommand\datapath{plots/qvi_membrane/data/Lshape}%
\includestandalone{plots/qvi_membrane/data/Lshape/qvi_membrane_est_final}%

\caption{Values of the error estimators for different values of $\gamma$. The refinement and update strategy (circles) is compared with the uniform refinement (squares).}\label{fig:qvi_membrane:Lshape:afem}
\end{figure}

\begin{figure}[h]
\newcommand\datapath{plots/qvi_membrane/data/slit}%
\includestandalone{plots/qvi_membrane/data/slit/qvi_membrane_est_final}%

\caption{Values of the error estimators for different values of $\gamma$. The adaptive refinement and update strategy (circles) is compared with the uniform refinement (squares).}\label{fig:qvi_membrane:slit:afem}
\end{figure}

\section{Conclusion}
In this work, we established a priori and a posteriori finite element techniques for a class of equilibrium problems \msa{and applied them}{} to quasi-variational inequalities. However, the developed framework is much more general and we expect the techniques to be well suitable to other \msa{instances}{} of equilibrium problems as well, such as \msa{e.g.}{} (generalized) Nash equilibrium problems.\\
In fact, many partial differential equations can be rewritten as first order conditions of a convex equilibrium problem and as such the derivation of error estimators \msa{may become more}{} straightforward \msa{or even feasible in first place, if other approaches fail}.

\begin{appendix}

\section{Appendix}

\begin{proof}[Proof of \refer{lem:my:conv}]
Let arbitrary sequences $\lambda_n \to \lambda$ and $u_n \rightharpoonup u$ be given \msa{and take}{} an arbitrary sequence $z_n \to z$ in $Z$\msa{.}{} \msa{Then,}{} for every $y \in Z$ there exists a sequence $y_n \to y$ such that
\begin{eq*}
\limsup_{n \to \infty} H(y_n,u_n) \leq H(y,u)
\end{eq*}
and we obtain
\begin{eq*}
\limsup_{n \to \infty} \my(\lambda_n,H(\ccdot,u_n))(z_n) &\leq \limsup_{n \to \infty} \left( H(y_n,u_n) + \frac{\lambda_n}{2}\|y_n - z_n\|^2 \right)\\
&\leq H(y,u) + \frac{\lambda}{2}\|y - z\|^2.
\end{eq*}
Since the choice of $y$ was arbitrary, we deduce
\begin{eq*}
\limsup_{n \to \infty} \my(\lambda_n,H(\ccdot,u_n))(z_n) \leq \my(\lambda,H(\ccdot,u))(z).
\end{eq*}
Let an arbitrary sequence $z_n \rightharpoonup z$ be given and take $y$ with $H(y,u) < \infty$ and $(y_n)_{\msa{n \in \N}} $ as above. Then, \msa{for every sequence of positive real numbers $(\eps_n)_{n\in \N}$ with $\eps_n \searrow 0$}{}, there exists a sequence $(\tilde y_n)_{n \in \N}$ such that
\begin{eq*}
\my(\lambda_n,\msa{H}(\ccdot,u_n))(z_n) &\leq H(\tilde y_n,u_n) + \frac{\lambda_n}{2}\|\tilde y_n - z_n\|_Z^2\\
&\leq \my(\lambda_n,\msa{H}(\ccdot,u_n))(z_n) + \eps_n.
\end{eq*}
Since $H$ is bounded from below and $\limsup_{n \to \infty} \my(\lambda_n,H(\ccdot,u_n))(z_n) \leq \my(\lambda, H(\ccdot,u))(z) < \infty$ and $\msa{(y_n)_{n \in \N}, (z_n)_{n\in\N}}$ are both bounded, we deduce the boundedness of \msa{$(\tilde y_n)_{n \in \N}$}. Hence, we can extract a subsequence \msa{$(n_k)_{k \in \N}$}{} with $\lim \my(\lambda_{n_k}, H(\ccdot, u_{n_k}))(z_n) = \liminf \my(\lambda_n, H(\ccdot, u_n))(z_n)$ such \msa{$(y_{n_k})_{k \in \N}$ is}{} weak\-ly convergent \msa{with}{} limit $y^*$ and obtain
\begin{eq*}
\my(\lambda, &H(\ccdot,u))(z) \leq H(y^*,u) + \frac{\lambda}{2}\|y^* - z\|^2_{\msa{Z}}\\
&\leq \liminf_{k \to \infty} \left( H(\tilde y_{n_k},u_{n_k}) + \frac{\lambda_{n_k}}{2}\|\tilde y_{n_k} - z_{n_k}\|^2_{\msa{Z}} \right)\\
&\leq \liminf_{k \to \infty} \left( \my(\lambda_{n_k}, H(\ccdot,u_{n_k}))(z_{n_k}) + \eps_{n_k} \right)\\
&= \liminf_{n \to \infty} \my(\lambda_n,H(\ccdot,u_n))(z).
\end{eq*}
From this, we deduce the assertion.
\end{proof}

\newcommand\Top{\ensuremath{\operatorname{Top}}}

\begin{lem}\label{lem:apx:derivative}
Consider convex, proper, lsc. functionals $f, \pi, \rho: X \to \R\cup \{\infty\}$,  with $f$ being subdifferentiable and $\dom \pi = \dom \rho = X$ \msa{and}{} $\pi, \rho$ being continuous \msa{as well as}{} $\rho$ being strongly convex. Define the functional
\begin{eq*}
v(\gamma) := \inf_{x \in X}\left(f(x) + \gamma \pi(x) + \frac{1}{\gamma}\rho(x)\right).
\end{eq*}
Then $v: (0,\infty) \to \R$ is differentiable \msa{and its first derivative reads}{}
\begin{eq*}
v'(\gamma) = \pi(x_\gamma) - \frac{1}{\gamma^2}\rho(x_\gamma)
\end{eq*}
with $x_\gamma := \argmin_{x \in X}\left( f(x) + \gamma \pi(x) + \frac{1}{\gamma} \rho(x) \right)$.
\end{lem}
\begin{proof}
We \msa{estimate}{} based on the definition of the value functional
\begin{eq*}
v(\gamma + \eta) &- v(\gamma) \leq f(x_\gamma) + (\gamma + \eta) \pi(x_\gamma) + \frac{1}{\gamma + \eta}\rho(x_\gamma)\\
&- f(x_\gamma) - \gamma \pi(x_\gamma) - \frac{1}{\gamma}\rho(x_\gamma)\\
&= \eta \pi(x_\gamma) - \frac{\eta}{\gamma(\gamma + \eta)} \rho(x_\gamma).
\end{eq*}
Dividing by $\eta > 0$ and passing to zero leads to
\begin{eq*}
\limsup_{\eta \searrow 0} \frac{v(\gamma + \eta) - v(\gamma)}{\eta} \leq \pi(x_\gamma) - \frac{1}{\gamma^2}\rho(x_\gamma)
\end{eq*}
and for $\eta < 0$ to
\begin{eq*}
\liminf_{\eta \nearrow 0} \frac{v(\gamma + \eta) - v(\gamma)}{\eta} \geq \pi(x_\gamma) - \frac{1}{\gamma^2} \rho(x_\gamma).
\end{eq*}
For the other \msa{direction}{} we use the first order systems reading as
\begin{eq*}
0 &\in \partial f(x_\gamma) + \gamma \partial \pi(x_\gamma) + \frac{1}{\gamma} \partial \rho(x_\gamma) &&\text{in } X^* \text{ and}\\
0 &\in \partial f(x_{\gamma + \eta}) + (\gamma + \eta) \partial \pi(x_{\gamma+\eta}) + \frac{1}{\gamma+\eta} \partial \rho(x_{\gamma+\eta}) &&\text{in } X^*.
\end{eq*}
\msa{Then, there exist}{} $\phi \in \partial f(x_\gamma), \phi' \in \partial f(x_{\gamma+\eta})$ as well as $p \in \partial \pi(x_\gamma), p' \in \partial \pi(x_{\gamma+\eta})$ and $r \in \partial \rho(x_\gamma), r' \in \partial \rho(x_{\gamma+\eta})$ \msa{which leads}{} to 
\begin{eq*}
0 &= \langle \phi' - \phi, x_{\gamma+\eta} - x_\gamma \rangle_{\msa{X^*,X}} + (\gamma+\eta) \langle p' - p, x_{\gamma+\eta} - x_\gamma \rangle_{\msa{X^*,X}{}}\\
&+ \frac{1}{\gamma+\eta}\langle r' - r, x_{\gamma+\eta} - x_\gamma \rangle_{\msa{X^*,X}} + \eta \langle p, x_{\gamma+\eta} - x_\gamma \rangle_{\msa{X^*,X}}\\
&+ \left( \frac{1}{\gamma+\eta} - \frac{1}{\gamma} \right)\langle r, x_{\gamma+\eta} - x_\gamma \rangle_{\msa{X^*,X}{}}
\end{eq*}
by testing with $x_{\gamma+\eta} - x_\gamma$. Using the strong convexity of $\rho$ we \msa{estimate}{}
\begin{eq*}
\frac{\alpha}{\gamma+\eta} \|x_{\gamma+\eta} - x_\gamma\|^2_{\msa{X}} \leq -\eta \langle p, x_{\gamma+\eta} - x_\gamma \rangle_{\msa{X^*,X}} - \frac{\eta}{\gamma(\gamma+\eta)}\langle r, x_{\gamma+\eta} - x_\gamma \rangle_{\msa{X^*,X}}
\end{eq*}
and thus
\begin{eq*}
\|x_{\gamma+\eta}-x_\gamma\|_{\msa{X}} \leq \frac{\eta}{\alpha} \left( (\gamma+\eta)\|p\|_{\msa{X^*}} + \frac{1}{\gamma} \|r\|_{\msa{X^*}} \right).
\end{eq*}
Hence, we can choose $\delta > 0$ with $|\eta| < \delta$ and make $\|x_{\gamma+\eta}-x_\gamma\|_{\msa{X}} < \eps$. Therefore, we obtain the continuity \msa{of}{} $\gamma \mapsto x_\gamma$ and subsequently
\begin{eq*}
v(\gamma+\eta) &- v(\gamma) \geq v(\gamma+\eta) - f(x_{\gamma+\eta}) - \gamma\pi(x_{\gamma+\eta}) - \frac{1}{\gamma}\rho(x_{\gamma+\eta})\\
&= \eta \pi(x_{\gamma+\eta}) - \frac{\eta}{\gamma(\gamma+\eta)} \rho(x_{\gamma+\eta}).
\end{eq*}
Using the previously proven continuity of the solution mapping we obtain for $\eta > 0$ that
\begin{eq*}
\liminf_{\eta \searrow 0} \frac{v(\gamma+\eta) - v(\gamma)}{\eta} \geq \pi(x_{\gamma}) - \frac{1}{\gamma^2} \rho(x_{\gamma})
\end{eq*}
and for $\eta < 0$ that
\begin{eq*}
\limsup_{\eta \nearrow 0} \frac{v(\gamma+\eta) - v(\gamma)}{\eta} \leq \pi(x_\gamma) - \frac{1}{\gamma^2} \rho(x_\gamma),
\end{eq*}
which yields the assertion.
\end{proof}

\begin{thm}[cf. \mbox{\cite[Theorem 2.1]{bib:BotConjugateDuality}}]\label{thm:fencheldual}
Let $X, Y$ be Banach spaces and let $f : X \to \bR$ and $g:Y \to \bR$ be proper and convex functionals and $A:X \to Y$ be a linear continuous operator. If the \emph{constraint qualification}
\begin{eq}[\tag{CQ}]\label{eq:fencheldual:cq}
0 \in \core{\dom g - A \dom f}
\end{eq}
holds, then also the identity
\begin{eq}[\tag{PD}]\label{eq:fencheldual:duality}
\inf_{x \in X} \left(\vphantom{\sum}f(x) + g(Ax)\right) + \min_{y^* \in Y^*} \left(\vphantom{\sum}f^*(-A^* y^*) + g^*(y^*)\right) = 0
\end{eq}
holds true.
\end{thm}

\begin{calc}\label{calc:conjugate}
Let $X$ be a reflexive Banach space and $f: X \to \bR$ be a \msa{proper}, convex, lower semi-continuous functional as described below. Then\msa{,}{} we calculate the conjugate functional $f^* : X^* \to \bR$:
\begin{enum}
\item Let $X = L^2(\Omega; \R^d)$ and $f(p):= \frac{1}{2}\|p\|_{L^2(\Omega;\R^d)}^2$, then $f^*(q) :=\\ \frac{1}{2}\|q\|_{L^2(\Omega;\R^d)}^2$.
\item Let $X$ be an arbitrary reflexive Banach space and define $f(x) := \langle \xi, x \rangle_{X^*,X}$ for some $\xi \in X^*$, then $f^*(x^*) := I_{\{\xi\}}(x^*)$.
\item Let $X = L^2(\Omega)$ and $f(y):= I_{K_{L^2}}(y)$ with $K_{L^2} := \{y \in L^2(\Omega) : y \geq \lowobs\}$ for some $\lowobs \in L^2(\Omega)$, then $f^*(z) = (z,\lowobs)_{\msa{L^2(\Omega)}} + I_{L^2_-(\Omega)}(z) $.
\item Let $X = L^2(\Omega)$ and $f(u):= I_{V_{L^2}}(u)$ with $V := \{z \in L^2(\Omega) : \int_\Omega (z - 1) \dx = 0 \}$, then $f^*(z) = I_\R(z) + \int_{\Omega} z \dx$.
\end{enum}
\end{calc}
\begin{proof}
\emph{ad (i):} See \cite[Example 13.6]{bib:BauschkeCombettes}.\\
\emph{ad (ii):} Take an arbitrary $x^* \in X^*$. If $x^* \neq \xi$, then there exists an element $x \in X$ such that $\langle x^* - \xi, x \rangle_{\msa{X^*,X}} > 0$ and hence we can take $x_n = n x$ to obtain $\lim_{n \to \infty} \left(\langle x^*, x_n \rangle_{\msa{X^*,X}} - f(x_n)\right) = \lim_{n \to \infty} \langle x^* - \xi, x_n \rangle_{\msa{X^*,X}} = \infty$, which implies $f^*(x^*) = \infty$. \msa{Otherwise,}{} $\langle x^* - \xi, x \rangle_{\msa{X^*,X}} = 0$ for all $x \in X$ and hence $f^*(\xi) = 0$, which proves the assertion.\\
\emph{ad (iii):} Take an arbitrary $z \in L^2(\Omega)$, then we obtain
\begin{eq*}
f^*(z) = \sup_{y \in L^2(\Omega)} \left( (z,y)_{\msa{L^2(\Omega)}} - I_{K_{L^2}}(y) \right) = (z, \lowobs)_{L^2(\Omega)} + \sup_{y \in L^2_+(\Omega)} (z,y)_{L^2(\Omega)}.
\end{eq*}
If $z \in L^2_-(\Omega)$, then the last term is zero. Otherwise, take the sequence $y_n = n z^+$ and we get $(z,y_n)_{\msa{L^2(\Omega)}} = n \int_\Omega z^{2+} \dx \to \infty$ as $n \to \infty$. Thus, we obtain $f^*(z) = (z,\lowobs)_{L^2(\Omega)} + I_{L^2_-(\Omega)}(z)$.\\
\emph{ad (iv):} Given $z \in L^2(\Omega)$ we can propose the decomposition $z = v + \alpha$ with $\alpha \msa{:=}{} \frac{1}{|\Omega|} \int_\Omega z \dx$ and $v \msa{:=}{} z - \alpha$. Then, we see $\alpha \in \R$ and $\int_\Omega v \dx = \int_\Omega z \dx - \int_\Omega \alpha \dx = \int_{\Omega} z \dx - \int_\Omega z\dx = 0$. Thus, we obtain
\begin{eq*}
f^*(z) &= \sup_{y \in L^2(\Omega)} \left( (z,y)_{\msa{L^2(\Omega)}} - I_{V_{L^2}}(y) \right) = \sup_{y \in V_{L^2}} \left( (v,y)_{\msa{L^2(\Omega)}} + \alpha \int_\Omega y \dx \right)\\
&= \sup_{y \in V_{L^2}} (v,y)_{\msa{L^2(\Omega)}} + \int_\Omega z \dx.
\end{eq*}
If now $v \neq 0$ we can take $y_n := n(v + \frac{1}{|\Omega|}) \in V_{L^2}$ and obtain $(v,y_n) = n\|v\|_{L^2(\Omega)}^2 + n \frac{1}{|\Omega|}\int_\Omega v \dx = n \|v\|^2_{\msa{L^2(\Omega)}} \to \infty$ as $n \to \infty$. Otherwise, we obtain zero for the first term. This is the same as $z = \frac{1}{|\Omega|} \int_\Omega z \dx$, respectively $z \in \R$. Hence, we deduce the assertion.
\end{proof}
\end{appendix}

\bibliographystyle{alpha}
\bibliography{mylibrary.bib}

\end{document}

%% file: plots/qvi_thermo/data/square/qvi_thermo_est_final.tex
\begin{tikzpicture}
\begin{axis}[xlabel=number of degrees of freedom,ylabel=eta2,xmode=log,ylabel=$\eta$,xmode=log,ymode=log,legend pos=outer north east,legend cell align={left}]
\addplot[color = red!25.0!purple, mark = square*] table {\datapath/est_data/est_uniform_0000.dat};
\addlegendentry{$\gamma =     1.0000$}
\addplot[color = red!50.0!purple, mark = square*] table {\datapath/est_data/est_uniform_0001.dat};
\addlegendentry{$\gamma =  1001.0000$}
\addplot[color = red!75.0!purple, mark = square*] table {\datapath/est_data/est_uniform_0002.dat};
\addlegendentry{$\gamma =  2153.0736$}
\addplot[color = red!100.0!purple, mark = square*] table {\datapath/est_data/est_uniform_0003.dat};
\addlegendentry{$\gamma =  5002.1582$}
\addplot[color = blue!14.285714285714285!green, mark = *] table {\datapath/est_data/est_adaptive_0000.dat};
\addlegendentry{$\gamma =     1.0000$}
\addplot[color = blue!28.57142857142857!green, mark = *] table {\datapath/est_data/est_adaptive_0001.dat};
\addlegendentry{$\gamma =  1001.0000$}
\addplot[color = blue!42.857142857142854!green, mark = *] table {\datapath/est_data/est_adaptive_0002.dat};
\addlegendentry{$\gamma =  2057.8030$}
\addplot[color = blue!57.14285714285714!green, mark = *] table {\datapath/est_data/est_adaptive_0003.dat};
\addlegendentry{$\gamma =  4701.9254$}
\addplot[color = blue!71.42857142857143!green, mark = *] table {\datapath/est_data/est_adaptive_0004.dat};
\addlegendentry{$\gamma = 12615.4908$}
\addplot[color = blue!85.71428571428571!green, mark = *] table {\datapath/est_data/est_adaptive_0005.dat};
\addlegendentry{$\gamma = 34919.8238$}
\addplot[color = blue!100.0!green, mark = *] table {\datapath/est_data/est_adaptive_0006.dat};
\addlegendentry{$\gamma = 101849.8974$}
\logLogSlopeTriangleDown{0.5}{-0.16}{0.2}{-0.5}{black}{$\frac{1}{2}$};
\logLogSlopeTriangleUp{0.85}{0.16}{0.4}{-0.666666666}{black}{$\frac{2}{3}$};
\end{axis}
\end{tikzpicture}

%% file: plots/qvi_membrane/data/Lshape/qvi_membrane_est_final.tex
\begin{tikzpicture}
\begin{axis}[xlabel=number of degrees of freedom,ylabel=eta2,xmode=log,ylabel=$\eta$,xmode=log,ymode=log,legend pos=outer north east,legend cell align={left}]
\addplot[color = red!25.0!purple, mark = square*] table {\datapath/est_data/est_uniform_0000.dat};
\addlegendentry{$\gamma =     1.0000$}
\addplot[color = red!50.0!purple, mark = square*] table {\datapath/est_data/est_uniform_0001.dat};
\addlegendentry{$\gamma =  1424.1338$}
\addplot[color = red!75.0!purple, mark = square*] table {\datapath/est_data/est_uniform_0002.dat};
\addlegendentry{$\gamma = 19361.7482$}
\addplot[color = red!100.0!purple, mark = square*] table {\datapath/est_data/est_uniform_0003.dat};
\addlegendentry{$\gamma = 38246.2827$}
\addplot[color = blue!16.666666666666664!green, mark = *] table {\datapath/est_data/est_adaptive_0000.dat};
\addlegendentry{$\gamma =     1.0000$}
\addplot[color = blue!33.33333333333333!green, mark = *] table {\datapath/est_data/est_adaptive_0001.dat};
\addlegendentry{$\gamma =  1414.7859$}
\addplot[color = blue!50.0!green, mark = *] table {\datapath/est_data/est_adaptive_0002.dat};
\addlegendentry{$\gamma = 20168.5364$}
\addplot[color = blue!66.66666666666666!green, mark = *] table {\datapath/est_data/est_adaptive_0003.dat};
\addlegendentry{$\gamma = 271096.9713$}
\addplot[color = blue!83.33333333333334!green, mark = *] table {\datapath/est_data/est_adaptive_0004.dat};
\addlegendentry{$\gamma = 324960.1151$}
\addplot[color = blue!100.0!green, mark = *] table {\datapath/est_data/est_adaptive_0005.dat};
\addlegendentry{$\gamma = 336633.0622$}
\logLogSlopeTriangleDown{0.6}{-0.16}{0.1}{-0.5}{black}{$\frac{1}{2}$};
\logLogSlopeTriangleUp{0.9}{0.16}{0.3}{-0.4}{black}{$\frac{2}{5}$};
\end{axis}
\end{tikzpicture}

%% file: plots/qvi_membrane/data/slit/qvi_membrane_est_final.tex
\begin{tikzpicture}
\begin{axis}[xlabel=number of degrees of freedom,ylabel=eta2,xmode=log,ylabel=$\eta$,xmode=log,ymode=log,legend pos=outer north east,legend cell align={left}]
\addplot[color = red!25.0!purple, mark = square*] table {\datapath/est_data/est_uniform_0000.dat};
\addlegendentry{$\gamma =     1.0000$}
\addplot[color = red!50.0!purple, mark = square*] table {\datapath/est_data/est_uniform_0001.dat};
\addlegendentry{$\gamma =  1001.0000$}
\addplot[color = red!75.0!purple, mark = square*] table {\datapath/est_data/est_uniform_0002.dat};
\addlegendentry{$\gamma =  9120.9029$}
\addplot[color = red!100.0!purple, mark = square*] table {\datapath/est_data/est_uniform_0003.dat};
\addlegendentry{$\gamma = 194342.0750$}
\addplot[color = blue!14.285714285714285!green, mark = *] table {\datapath/est_data/est_adaptive_0000.dat};
\addlegendentry{$\gamma =     1.0000$}
\addplot[color = blue!28.57142857142857!green, mark = *] table {\datapath/est_data/est_adaptive_0001.dat};
\addlegendentry{$\gamma =  1001.0000$}
\addplot[color = blue!42.857142857142854!green, mark = *] table {\datapath/est_data/est_adaptive_0002.dat};
\addlegendentry{$\gamma =  2569.0129$}
\addplot[color = blue!57.14285714285714!green, mark = *] table {\datapath/est_data/est_adaptive_0003.dat};
\addlegendentry{$\gamma =  7472.2564$}
\addplot[color = blue!71.42857142857143!green, mark = *] table {\datapath/est_data/est_adaptive_0004.dat};
\addlegendentry{$\gamma = 33201.3822$}
\addplot[color = blue!85.71428571428571!green, mark = *] table {\datapath/est_data/est_adaptive_0005.dat};
\addlegendentry{$\gamma = 66514.7056$}
\addplot[color = blue!100.0!green, mark = *] table {\datapath/est_data/est_adaptive_0006.dat};
\addlegendentry{$\gamma = 73839.2624$}
\logLogSlopeTriangleDown{0.6}{-0.16}{0.1}{-0.5}{black}{$\frac{1}{2}$};
\logLogSlopeTriangleUp{0.85}{0.16}{0.47}{-0.25}{black}{$\frac{1}{4}$};
\end{axis}
\end{tikzpicture}